\documentclass[12pt]{article}

\usepackage{amsmath, amsthm, amsfonts}

\usepackage[utf8]{inputenc}

\usepackage{graphicx}

\usepackage[T1]{fontenc}
\newtheorem*{theorem}{Theorem}
\newtheorem*{axiom}{Axiom}
\newtheorem*{corollary}{Corollary}
\newtheorem{remark}{Remark}

\title{An Elementary System of Axioms for Euclidean Geometry \\ based on Symmetry Principles}
\date{}
\author{Boris Čulina\\
  \small{Department of Mathematics, University of Applied Sciences Velika Gorica},\\
  	 \small{Zagreba\v{c}ka cesta 5, Velika Gorica, CROATIA}\\
 \small{email: boris.culina@vvg.hr}}

\begin{document}

\maketitle
\newpage
\textbf{Abstract.} In this article I develop an elementary system of axioms for Euclidean geometry. On  one hand, the system is based on the symmetry principles which express our a priori ignorant approach to space: all places are the same to us (the homogeneity of space), all directions are the same to us (the isotropy of space) and all units of length we use to create geometric figures are the same to us (the scale invariance of space).  On the other hand,  through the process of  algebraic simplification, this system of axioms directly provides the Weyl's system of axioms for Euclidean geometry. The system of axioms, together with its a priori interpretation, offers new views to philosophy and pedagogy of mathematics: (i) it supports the thesis that Euclidean geometry is a priori, (ii) it supports the thesis that in modern mathematics the Weyl's system of axioms is  dominant   to the Euclid's system  because it reflects the  a priori underlying symmetries, (iii) it gives a new and promising approach to learn geometry which, through the Weyl's system of axioms, leads  from the essential geometric symmetry principles of the mathematical nature directly to  modern mathematics.

\textbf{keywords:} symmetry, Euclidean geometry, axioms, Weyl's axioms, philosophy of geometry, pedagogy of geometry

\section{Introduction}
\label{intro}
The connection of Euclidean geometry with symmetries has a long history. I will be concerned here with the three symmetry principles: (i) the homogeneity of space: all places are the same,  (ii) the isotropy of space: all directions are the same, (iii) the scale invariance of space: all units of length we use to create geometric figures are the same. In 17th century John Wallis  proved, assuming other Euclid's postulates, that  the scale invariance principle "For every figure there exists similar figure of arbitrary magnitude."  is equivalent to the Euclid's fifth postulate \cite{Wallis}. Wallis considered  his postulate to be more convincing than Euclid's fifth postulate. Tracing back to the famous Riemann lecture at G\"{o}ttingen in 1854 \"{U}ber die Hypothesen welche der Geometrie zu Grunde liegen" (\cite{Riemann}), it is well known that among all Riemann manifolds Euclidean geometry is characterized by the three symmetry principles. However, this characterisation is not an elementary one  because it presupposes the whole machinery of Riemann manifolds. The descriptions based on Klein's program of characterizing geometries by their principal groups of transformations (\cite{Klein}) also have non elementary character. As I am aware, there is no an elementary description (a description in terms of intuitive relations between points) of Euclidean geometry that is based on the three symmetry principles. Here I develop a system of axioms that provides such an elementary description. 

The importance and validity of the three symmetry principles has been recognized a long time ago. I offer a new interpretation of these principles. William Kingdon Clifford in \cite{Clifford1} and \cite{Clifford2} considers the three symmetry principles as the most  essential geometrical assumptions. He considers  that the principles are based on observations of the real space. Hermann von Helmholtz has the same opinion for the first two symmetry principles which he unifies in his principle of the free mobility of rigid bodies (\cite{Helmholtz}). Henri Poincar\'{e}, in his analysis of the real space \cite{Poincare},  comes to the conclusion that the first two symmetry principles are the most essential properties of the so called geometric space which  for him is not the real space but a ''conventional space'' -- the most convenient description of the real space. An interesting explanation of the validity of the three symmetry principles comes  from Joseph Delboeuf (\cite{Delboeuf}). He considers what remains when  we ignore all differences of things caused by their movements and mutual interactions. According to Delboeuf, in the ultimate abstraction from all diversities of real things we gain the homogeneous, isotropic, and scale invariant  space - the true geometric space which is Euclidean and which is different from the real space. My view on the principles is a different one: they are not a posteriori, the result of analysing the real space, but they are a priori. However, in my interpretation they are not a priori in Kant's sense, as an a priori form of rational cognition, but in a sense of our a priori ignorant approach to space: all places are the same to us (the homogeneity of space), all directions are the same to us (the isotropy of space) and all units of length we use for constructions in space are the same to us (the scale invariance of space). Because the system of axioms I develop here is based on the three symmetry principles and because it is equivalent to other systems of axioms for Euclidean geometry, under my a priori interpretation of the principles, the system  supports the thesis that Euclidean geometry is a priori, in the same way as number systems are a priori, the result  of modelling, not the world, but our activities in the world. This conclusion could satisfy Gauss who expressed his dissatisfaction with the epistemic status of Euclidean geometry in a letter to Olbers (\cite{Gauss}): ''I am ever more convinced that the necessity of our geometry cannot be proved, at least not by, and not for, our human understanding. Maybe in another life we shall attain insights into the essence of space which are now beyond our reach. Until then we should class geometry not with arithmetic, which stands purely a priori, but, say, with mechanics''.

The system of axioms I develop here directly provides, through the process of  algebraic simplification,  the Weyl's system of axioms for Euclidean geometry (\cite{Weyl})\footnote{brief description of the Weyl's system is on pages \pageref{Wa} and \pageref{We} }. Thus, the system of axioms I develop here, together with my a priori interpretation of the system, offers an explanation of the fact that  in modern mathematics the Weyl's system of axioms    is  dominant   to the Euclid's system. Although  the Euclid's system of axioms for Euclidean geometry is thought in school, the Weyl's system of axioms is used in modern mathematics, physics and engineering. Only through  the Weyl's system of axioms do we find the Euclidean structure in complex mathematical structures, and this enables us  to make our reasoning about them more visual and efficient. Today,  the Weyl's system of axioms is one of the essential  synthesizing tools of  modern mathematics while Euclid's system is of a secondary importance. My thesis is that this happens because  the original Euclidean system of axioms reflects a posteriori   intuition (hence physical intuition) about Euclidean geometry while  the Weyl's system of axioms reflects   a priori  intuition (hence mathematical intuition) about the underlying symmetries.

The system of axioms I develop here also gives a new and promising approach to learn Euclidean geometry  because (i) it gives an elementary and essential description of Euclidean geometry, (ii) it gives, through the process of  algebraic simplification, the  Weyl's system of axioms for Euclidean geometry which is  essential for  modern mathematics. Therefore, it opens the possibility to learn geometry  in a way which  leads  from the essential geometric symmetry principles (and these are of the mathematical nature) directly to  modern mathematics.

I gradually introduce a system of axioms about points labelled with letter ''A'', which have an immediate justification in intuitive ideas about relations between points and in the three symmetry principles supported by an idea of continuity of space. Then, I deduce from them  The Weyl's system of axioms about points and vectors labelled with letter ''W'', which are indeed equivalent to A -- axioms.  The complete list of axioms is displayed in the Appendix.

The primitive terms of the system of axioms are: (i)  equivalence of pairs of points (arrows), (ii) multiplication of a pair of points by a real number and (iii ) distance between points. The multiplication could be avoided.  However, the procedure to define the multiplication is somewhat lengthy and I prefer to introduce the multiplication as a new primitive term. Also, it is more simple to introduce the distance function (to add an arbitrary unit of measurement) as a new primitive  than to introduce  congruence between pairs of points as a new primitive term and define the distance function relative to the choice of a unit of measurement.

\section{Equivalence of pairs of points}
\label{sec:1}
Geometrical space $S$ will be modelled as a non empty set of objects termed \textbf{points}. The basic geometrical relation is the  position of one point relative to another (not necessarily different) point. That the position of a point $B$ relative to a point $A$ is the same as the position of a point $B'$ relative to a point $A'$ we will denote $AB\sim A'B'$ and we will say that  pairs or arrows $AB$ and $A'B'$ are \textbf{equivalent}. This is the first  primitive term of our system. It expresses a basic intuitive idea about relation of two points. The idea itself \textit{to be in the same relative position} implies that it is a relation of equivalence. This is the content of the first axiom.

\begin{axiom}[\textbf{A1}] $\sim$ is an equivalence relation. \end{axiom}

\noindent In more detail, it means:

\begin{axiom}[\textbf{A1.1}]  $AB\sim AB$  \hspace{5mm} (reflexivity)\end{axiom}

\begin{axiom}[\textbf{A1.2}]  $AB\sim A'B' \ \rightarrow \ A'B'\sim AB$  \hspace{5mm} (symmetry)\end{axiom}

\begin{axiom}[\textbf{A1.3}]   $AB\sim A'B' \ \land\  A'B'\sim A''B''  \ \rightarrow \ AB\sim A''B''$  \hspace{5mm} (transitivity)\end{axiom}

\noindent Because of an elementary character of these properties of $\sim$,  usually I will not mention them  in proofs of theorems.

Concerning a fixed point $A$ we can easily describe the equivalance relation $\sim$: by the very idea of the relative position of points, different points have different relative positions to $A$:

\begin{axiom}[\textbf{A2}]  $AB\sim AC  \ \rightarrow \  B=C$. \end{axiom}

Fundamental operations with arrows are to invert an arrow and to add an arrow to another arrow. The definitions follow:

\vspace{2mm}
\noindent \textbf{inverting arrow}:\hspace{5mm} $AB \ \mapsto \  -AB=BA$

\vspace{2mm}
\noindent \textbf{addition of arrows};\hspace{5mm} $AB, \  BC  \ \mapsto \  AB+BC=AC$
\vspace{2mm}

Because of axiom A2 we can extend addition of arrows:

\vspace{2mm}
\noindent \textbf{generalized addition of arrows};\hspace{5mm} $AB+CD=AB+BX$, where $BX\sim CD$, under the condition that there is such a point $X$.
\vspace{2mm}

\noindent By the homogeneity principle, the operations are invariant under the equivalence of arrows:

\begin{axiom}[\textbf{A3.1}]  $AB\sim A'B'  \ \rightarrow \  BA\sim B'A'$. \end{axiom}

\begin{axiom}[\textbf{A3.2}]  $AB\sim A'B' \ \land \ BC\sim B'C'  \ \rightarrow \  AC\sim A'C'$ \end{axiom}

Until now, we know only that $AB$ is equivalent to itself (reflexivity of $\sim$) and to no other arrow from the point $A$ (axiom A2). All other axioms are conditional statements. It remains to describe the equivalence of arrows originating from different points. 

\section{Multiplication of an arrow by a number}

\noindent In the next section I will sketch how  can we define the multiplication of an arrow by a number and establish  basic properties of the operation.  Although such  approach is conceptually more satisfactory, it is technically   very lengthy. Instead,   I will  directly postulate the properties of multiplication of arrows by real numbers, without losing a clear basis of the three symmetry principles. I introduce \textbf{multiplication of an arrow by a real number} as a new primitive operation based on an idea of stretching arrows and of an idea of iterative addition of the same arrow (numbers will be labelled with letters from the Greek alphabet):

\vspace{2mm}

$\cdot : \mathbb{R}\times S^2 \rightarrow S^2$ \hspace{1cm} $\lambda , A, B \ \mapsto \ \lambda \cdot AB$

\vspace{2mm}

\noindent Sometimes, since it is a common convention, we will not write the multiplication sign at all.

The very idea of the multiplication as stretching arrows is formulated in the next axiom:

\begin{axiom}[\textbf{A4}]  $\forall \lambda , A, B \ \exists C  \ \ \lambda \cdot AB=AC$. \end{axiom}

By the homogeneity principle,  multiplication of an arrow by a number is invariant under the equivalence of arrows:

\begin{axiom}[\textbf{A5}]  $AB\sim CD \ \rightarrow \ \lambda AB\sim \lambda CD$. \end{axiom}

For a  point $C$ such that   $AC = \lambda \cdot AB$ we will say that it is \textbf{along} $AB$. Also, for arrow $A$  we will  say that it is \textbf{along} $AB$. 

The very idea of the multiplication as addition of the same arrow  leads to the next axiom

\begin{axiom}[\textbf{A6.1}]  $1\cdot AB = AB$. \end{axiom}

By the homogeneity principle, we can translate any arrow along $AB$ to any point along $AB$. So, we can add such arrows. Specially,  we can add $\lambda\cdot AB$ and $\mu\cdot AB$ and the result will be $\lambda\cdot AB +\mu\cdot AB =\nu \cdot AB$ for some number $\nu$. Moreover, by the very idea of the multiplication as  iterative addition of the same arrow $\nu = \lambda + \mu$. This is the content of the next axiom:

\begin{axiom}[\textbf{A6.2}]  $ \lambda\cdot AB +\mu\cdot AB = (\lambda +\mu)\cdot AB$.\end{axiom}

\noindent Let's note that with this equation we  postulate also that the left side of the equation is defined. 

If we  stretch an arrow along $AB$ the result will be an arrow along $AB$, too. So, $\lambda\cdot(\mu\cdot AB)= \nu\cdot AB$, for some number $\nu$.Moreover, by the very idea of the multiplication as  iterative addition of the same (stretched) arrow $\nu = \lambda \cdot \mu$. This is the content of the next axiom: 

\begin{axiom}[\textbf{A6.3}]  $\lambda\cdot (\mu\cdot AB)=(\lambda\cdot\mu)\cdot AB$. \end{axiom}

\noindent Let's note that with this equation we  postulate also that $\lambda\cdot (\mu\cdot AB)$ is along $AB$.

Later, we will need to translate an arrow on its both sides:

\begin{corollary}[\textbf{A'4}]  $\forall A,B \ \exists D \ \ AB \sim BD$. \end{corollary}

\begin{proof}
	$2\cdot AB = (1+1)\cdot AB = \textrm{(by A6.2)\ } 1\cdot AB + 1\cdot AB = \textrm{(by A6.1)\ } AB + AB$. By the meaning of addition there is the point $D$ such that $AB\sim BD$. 
\end{proof}

\begin{corollary}[\textbf{C1}]  $\forall A,B \ \exists C \ \ CA\sim AB $. \end{corollary}
\begin{proof}
	Applying A'4 on $BA$  we can find  the point $C$ such that $BA\sim AC$. Because inverting arrows is invariant under $\sim$ (A3.1) it follows that $AB\sim CA$.	  
\end{proof}

The last axiom expresses  the scale invariance principle.  

\begin{figure}[h]
\begin{center}
\includegraphics[height=4cm]{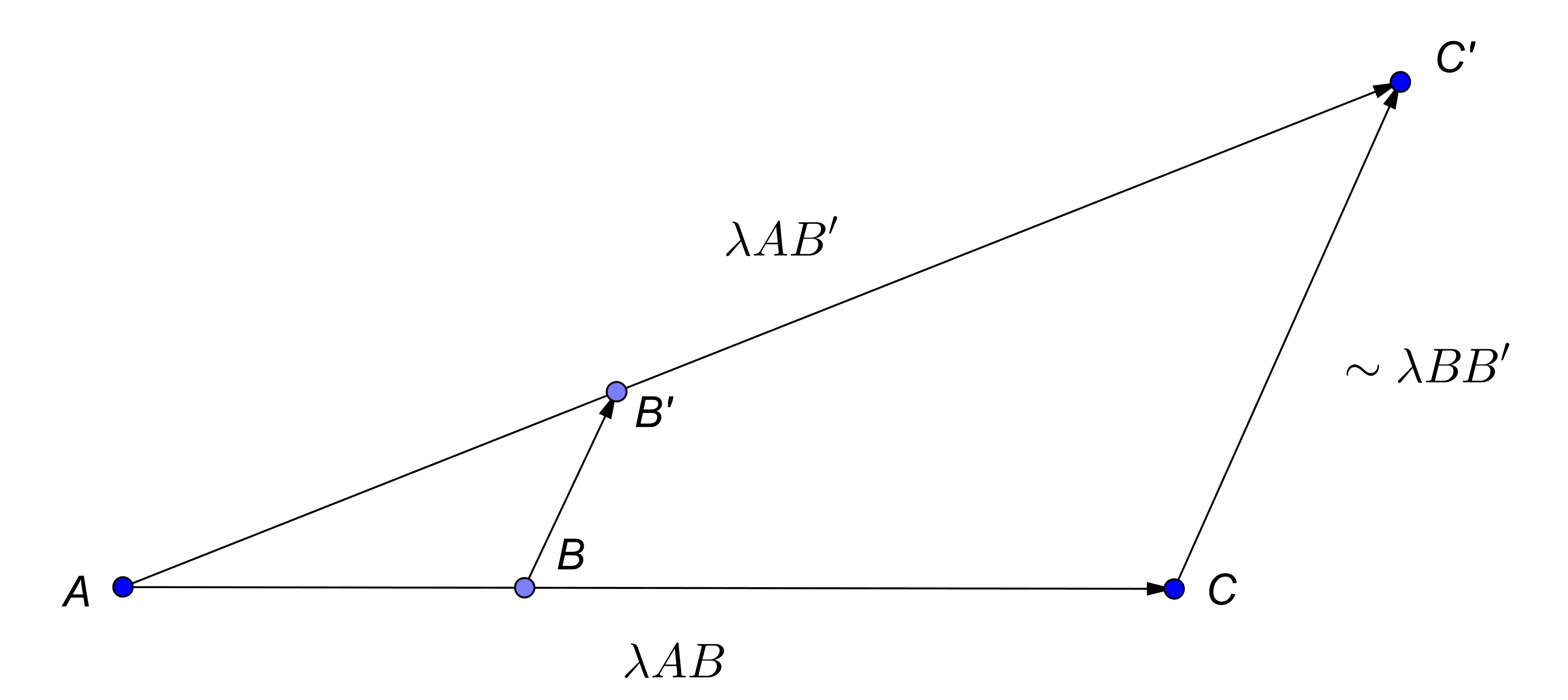}
\caption{ }
\label{fig:0}
\end{center}
\end{figure}

\begin{axiom}[\textbf{A7}] (the scale invariance axiom) If $ AC = \lambda\cdot AB$ and $ AC' = \lambda\cdot AB'$ then $CC' \sim \lambda\cdot BB'$. \normalfont{(Fig.\ref{fig:0})}   \end{axiom}


\begin{corollary}[\textbf{C2}]  $\lambda\cdot (AB+BB')=\lambda\cdot AB+\lambda\cdot BB'$. \end{corollary}

\begin{proof}
	By A4 there are points $C$ and $C'$ such that  $\lambda\cdot AB = AC$ and  $\lambda\cdot AB' = AC'$. Then, by axiom A7, $CC'\sim \lambda\cdot BB'$. Now, we calculate:
	$\lambda\cdot AB+\lambda\cdot BB' = AC + CC' = AC' = \lambda\cdot AB' = \lambda\cdot (AB+BB')$.
\end{proof}

\begin{figure}[h]
	\begin{center}
		\includegraphics[height=4cm]{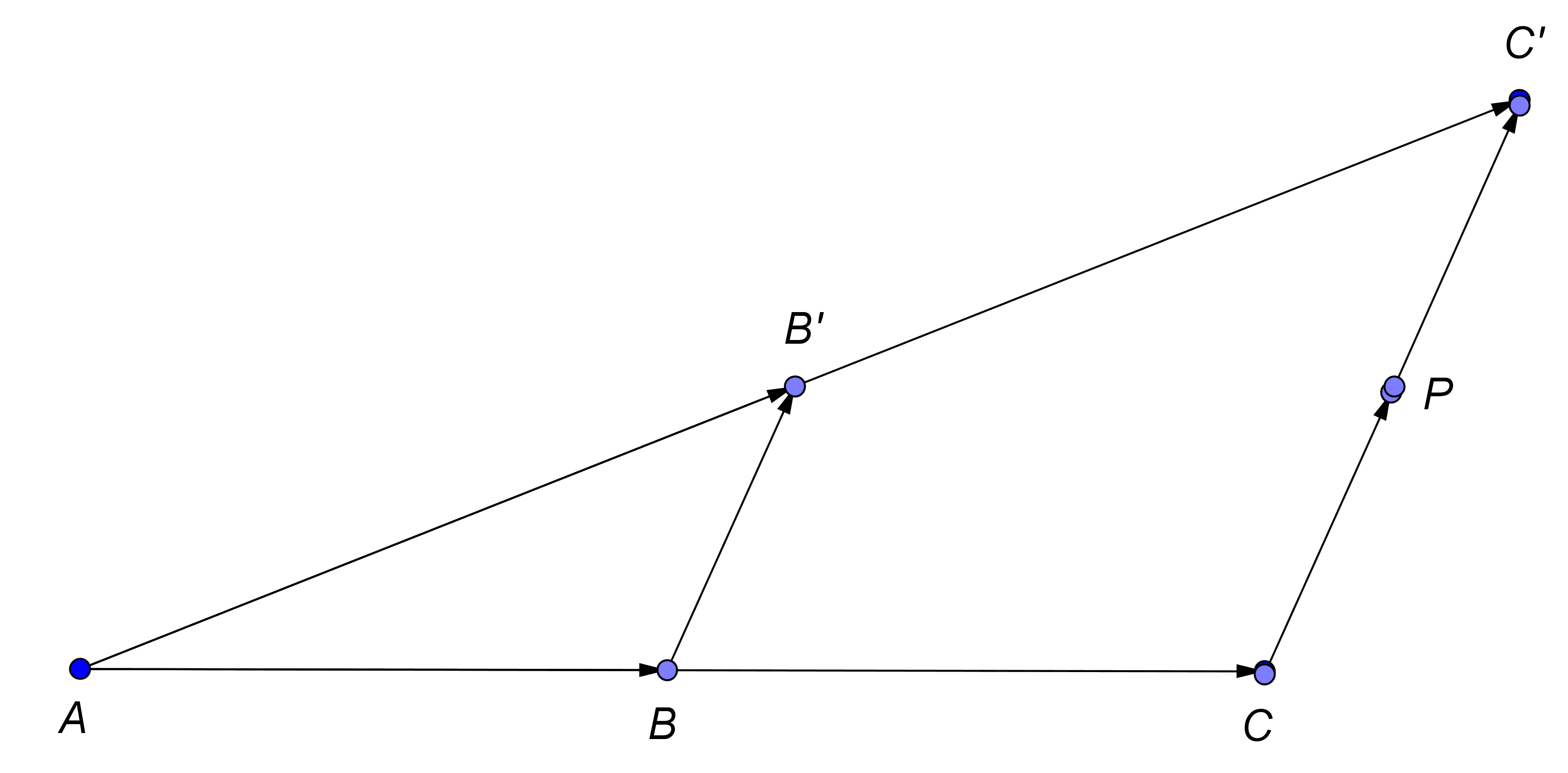}
		\caption{ }
		\label{fig:00}
	\end{center}
\end{figure}

Of the special interest is  a somewhat modified special case of the scale invariance axiom, for $\lambda = 2$:

\begin{theorem}[\textbf{A'5}]  (the elementary scale invariance law) $AB\sim BC$  and $AB'\sim B'C'   \ \rightarrow \ \exists \ P \ \ CP\sim PC' \sim BB'$. \normalfont{(Fig.\ref{fig:00})} \end{theorem} 

\begin{proof}
	Suppose that  $AB\sim BC$  and $AB'\sim B'C'$. Then  $ AC = AB+BC \sim  AB+AB = \textrm{(by A6.1) \ } 1\cdot AB+1\cdot AB =  \textrm{(by A6.2) \ }  2\cdot AB$. In the same way we prove that $ AC' = 2\cdot AB'$. Then, by A7, $CC' \sim 2\cdot BB'$. Multiplying this relation with $\frac{1}{2}$, by A5 and A6.3, it follows that $BB'\sim \frac{1}{2}\cdot CC'$. By  A4, there is the point $P$ such that $\frac{1}{2}\cdot CC' = CP$. Therefore,  $CP\sim PC'\sim BB'$.
\end{proof}

This theorem has two important consequences.

\begin{theorem}[\textbf{T3}] (the  unique translation of arrows law) $\forall A,B,A'\ \exists ! B' \ \ AB\sim A'B'$. \end{theorem}

\begin{proof}
	We will accomodate notation to the elementary scale invariance law. We will prove that for all $B, B', C$ there exists the unique point P such that $BB'\sim CP$. By C1 there is the point $A$ such that $AB\sim BC$. By A'4 there is the point C' such that $AB'\sim B'C'$. By  A'5 there is the point P such that $BB'\sim CP$. The point $P$ is a unique such point. If there were another point $P'$ such that $BB'\sim CP'$, then by symmetry (A1.2) and transitivity (A1.3) of $\sim$, we infer $CP'\sim CP$. From this statement, by A2, it follows that $P'=P$.
\end{proof}

The  unique translation of arrows law enables us to add arbitrary arrows, without any condition, as we have done before.

\vspace{2mm}
$AB+CD=AB+BX$, where $BX\sim CD$
\vspace{2mm}

\begin{figure}[h]
	\begin{center}
		\includegraphics[height=3.5cm]{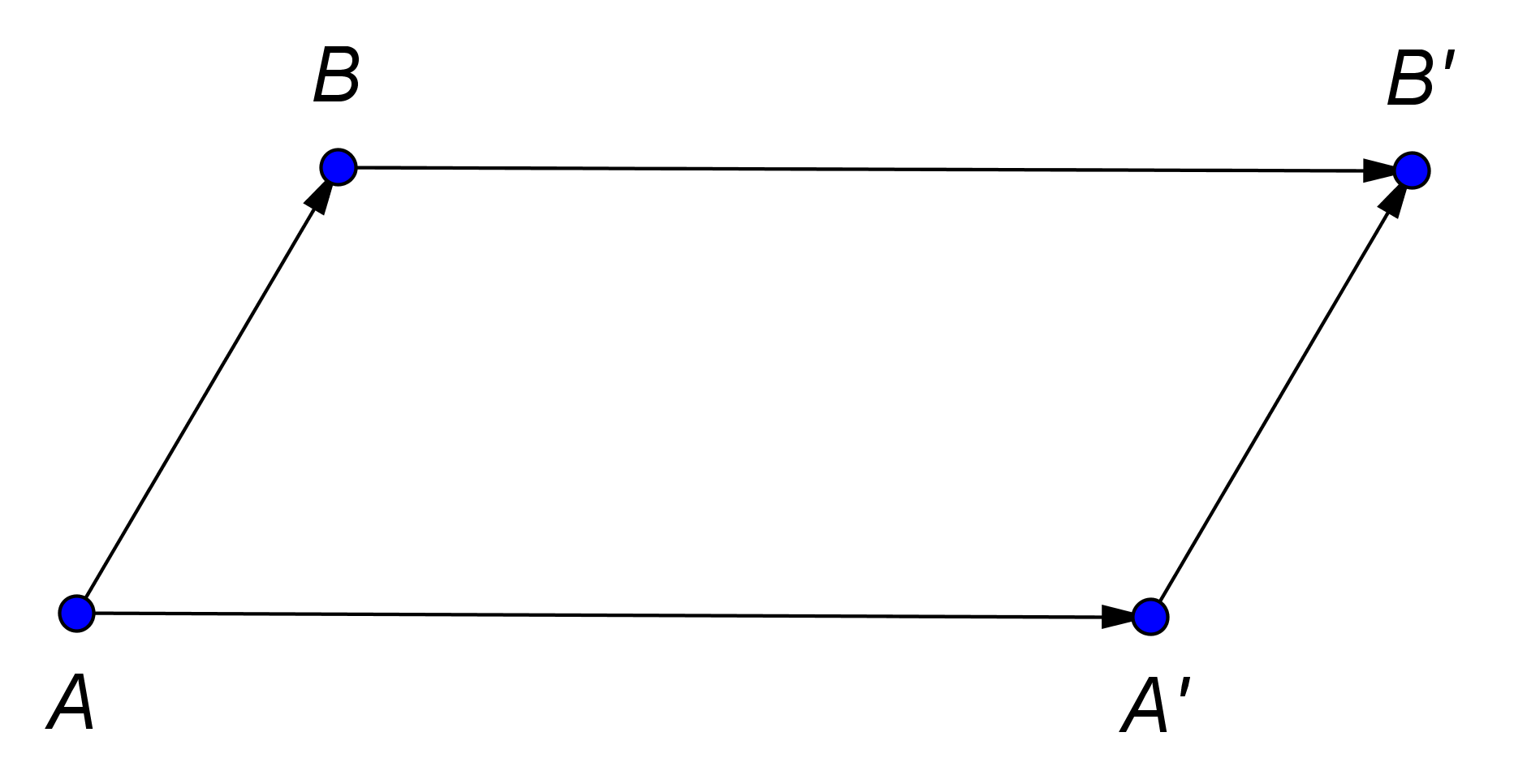}
		\caption{ }
		\label{fig:2}
	\end{center}
\end{figure}

\begin{theorem}[\textbf{T4}] (the parallelogram law) $AB\sim A'B'   \ \rightarrow \ AA'\sim BB'$. \normalfont{(Fig.\ref{fig:2})}\end{theorem}

\begin{proof}
	We will accomodate notation to the elementary scale invariance axiom. We will prove that if $BB'\sim CP$ then $BC\sim B'P$. By C1 there is the point $A$ such that $AB\sim BC$. By A'4 there is the point C' such that $AB'\sim B'C'$. By  A'5 there is the point P' such that 	$CP'\sim P'C' \sim BB'$. By transitivity of $\sim$, we infer $CP'\sim CP$. So, by A2, $P=P'$. Therefore, $BB'\sim PC$. Now, we calculate:
	
	$BC \sim AB = AB'+B'B \sim \ \textrm{(by A3)} \ B'C'+C'P \sim B'P$
\end{proof}	

\begin{remark}
	Assuming only axioms A1 to A3 and corollary A'4,  without any axiom on the multiplication, we can prove that the elementary scale invariance axiom is equivalent to the conjunction of the  unique translation of arrows law and parallelogram law.
\end{remark}

\begin{remark}
	Until now we haven't used axiom A1.1 (reflexivity). Indeed we can prove it: by the unique translation of arrows law T3, for given points $A$ and $B$ there is the unique $C$ such that $AB\sim AC$. From A2 it follows $B=C$. Therefore,  $AB\sim AB$.
\end{remark}
\begin{remark}
	The parallelogram law is very efficient algebraically in proofs. For example, assuming only the parallelogram law,  we can prove the equivalence of:  (i) A1.1 (reflexivity of  $\sim$) and the next corrolary C3 (all arrows of the type  $AA$ are equivalent), (ii) A1.2 (symmetry of $\sim$) and A3.1 (invariance under  $\sim$ of inverting arrows) ,  (iii)  A1.3 (transitivity of  $\sim$) and A3.2 (invariance under  $\sim$ of adding arrows).
\end{remark}

\begin{corollary}[\textbf{C5}]  $\forall A,B \ \ AA\sim BB$. \end{corollary}

\begin{proof}
	From reflexivity of $\sim$, we infer $AB\sim AB$. By  the parallelogram law,  it follows that $AA\sim BB$.
\end{proof}

\section{Elimination of the multiplication of an arrow by a number}

Here, I will sketch how  can we define the multiplication of an arrow by a number and establish  basic properties of the operation. I continue exposition from the second section, based on axioms A1 to A3, neglecting the third section. In the second section we have described the relation $\sim$ between arrows originating from the same point. How can we compare arrows originating from different points? There is an obvious case: we can translate an arrow $AB$ along itself.  In that way, by the homogeneity principle, we can reproduce the  position of $B$ relative to $A$ in the point $B$. This is the content of the axiom A'4 (which we have proved in the second section from the multiplication axioms). 

\begin{axiom}[\textbf{A'4}]  $\forall A, B \ \exists   C \ \  AB\sim BC$. \end{axiom}

An immediate consequence is that we can translate $AB$ in the opposite direction:

\begin{corollary} [\textbf{C1}]
	$\forall A, B \ \exists   D \ \  DA\sim AB  $
\end{corollary}.

\noindent We have proved this in the previous section using only A'4, A3.2 and A1.2.

As the next axiom we will use the elementary scale invariance axiom (which we have also proved in the previous section from the multiplication axioms)

\begin{axiom}[\textbf{A'5}]  (the elementary scale invariance axiom) $AB\sim BC$  and $AB'\sim B'C'   \ \rightarrow \ \exists \ P \ \ CP\sim PC' \sim BB'$. \end{axiom}


These two axioms and axioms A1 to A3 yield   the  unique translation of arrows law T3 and the parallelogram law T4. I have used nothing else in the proofs of these laws in the previous section.

Because  the  unique translation of arrows law  enables us to add arbitrary arrows, we can add  an arrow to itself, as a special case. It enables us to multiply the arrow by a positive natural number:

\vspace{2mm}
$1\cdot AB=AB, \ 2\cdot AB=AB+AB, \ 3\cdot AB=2\cdot AB+AB,\ldots$
\vspace{2mm}

\noindent In that way we can define recursively the multiplication of an arrow by an arbitrary positive natural number.   Then, we can define the multiplication of an arrow by an arbitrary integer:

\vspace{2mm}
$0\cdot AB=AA$,
\vspace{2mm}
$(-n)\cdot AB =n\cdot (-AB)$, where $n$ is a positive natural number.
\vspace{2mm}
\noindent Furthermore, we can introduce axiomatically  (and justify it by the symmetry principles and an idea of divisibility of space) the existence of the midpoint  $P$ of an arrow $AB$, as the point for which $AP\sim PB$,

\begin{axiom}[\textbf{A'6}]  $\forall A, B \ \exists  ! P \ \  AP\sim PB$. \end{axiom}

\noindent and then define multiplication of the arrow by $\frac{1}{2}$:

\vspace{2mm}
$\frac{1}{2}\cdot AB =AP$
\vspace{2mm}

\noindent By repeating bisection we can define recursively multiplication of an arrow by $\frac{1}{2^n}$. On that basis we can define multiplication of an arrow by an arbitrary rational number of the form $\frac{m}{2^n}$, where $m$ is an integer and $n$ is a positive natural number.  

We can extend the multiplication to all real numbers using an appropriate axiom of continuity. For example, 

\begin{axiom}[\textbf{A'7}] Let $p_2 (A,B)$, for $A\neq B$, be  the set of all points C such that $AC = \lambda \cdot AB$, where $\lambda$ is rational number of the form $\frac{m}{2^n}$. Let's define the order relation $<_2$ on $p_2 (A,B)$: for   $AC = \lambda _C \cdot AB$ and   $AD = \lambda _D \cdot AB$, $C<D \ \leftrightarrow \  \lambda _C < \lambda _D$. Then there is the unique linearly ordered extension $p (A,B)$  of $p_2 (A,B)$ such that  every  bounded above non empty subset of $p (A,B)$  has supremum.\end{axiom}

\noindent  Now, we can  define multiplication of an arrow by an arbitrary real number $\lambda$. Because for every $\lambda >0$ there exists the unique strictly increasing sequence $\frac{m_n}{2^n}$ such that $\lambda = \lim_{n \to \infty} \frac{m_n}{2^n}$, we define

\vspace{2mm}
$\lambda \cdot AB = A\sup\{ C_n | AC_n = \frac{m_n}{2^n}\cdot AB\}$

$(-\lambda) \cdot  AB =-(\lambda) \cdot  AB  $
\vspace{2mm}

\noindent By the homogeneity principle, all these operations are invariant under $\sim$:

\begin{axiom}[\textbf{A'8}]  The operations of bisection and supremum are invariant under $\sim$. \end{axiom}

From these axiom all the axioms for the multiplication can be proved.

\section{Vectors}

Axiom A1 enables us to define vectors. Since,  by A1,  $\sim$ is an equivalence relation, it classifies arrows (pairs of points) into  classes of mutually equivalent arrows. We define \textbf{vectors} as these  equivalence classes. The set of all vectors will be denoted  $\overrightarrow{S}$. To every pair of points $AB$ we will associate the vector $\overrightarrow{AB}$, the equivalence class  to which $AB$ belongs:

\vspace{2mm}
$\overrightarrow{AB}=\{CD | CD\sim AB\}$
\vspace{2mm}

\noindent So, $\overrightarrow{\ }$ maps pairs of points to vectors: $\ ^{\overrightarrow{ \hspace{4mm} }}: S^2 \rightarrow \overrightarrow{S}$.

For every pair of points (arrows) from a vector (an equivalence class)  we  say that it \textbf{represents} the vector. Thus, for example, a pair $AB$ represents the vector $\overrightarrow{AB}$.

\begin{theorem}[\textbf{W1}] For every point  $A$  the function $X \ \mapsto \ \overrightarrow{AX}$ is a bijection from  set of points $S$ onto  set of vectors $\overrightarrow{S}$. \end{theorem}

\begin{proof} The claim is an immediate consequence of  the unique translation of arrows law (T3).\end{proof}

For a given point $A$ the inverse function of  bijection  $X \ \mapsto \ \overrightarrow{AX}=\overrightarrow{x}$ maps every vector  $\overrightarrow{x}$ to the point $X$ which we will denote  $X=A+\overrightarrow{x}$. Thus,

\vspace{2mm}
$A+\overrightarrow{x}=X \ \leftrightarrow \ \overrightarrow{AX}=\overrightarrow{x}$
\vspace{2mm}

We will transfer operations with arrows into operations with corresponding vectors, in a way invariant under  relation $\sim$.

\vspace{2mm}
\noindent \textbf{null vector}:\hspace{5mm} $\overrightarrow{0}=\overrightarrow{AA}$
\vspace{2mm}

\noindent By  corollary C5, the definition is correct because it does not depend on the choice of a point  $A$.

\vspace{2mm}
\noindent \textbf{inverse vector}:\hspace{5mm} $-\overrightarrow{AB}=\overrightarrow{BA}$
\vspace{2mm}

\noindent By  axiom A3.1 (invariance under $\sim$ of inverting arrows), the definition is correct because it does not depend on the choice of an arrow $AB$.

\vspace{2mm}
\noindent \textbf{Addition of vectors} 
(\textbf{W2}):  \hspace{5mm} $\overrightarrow{AB}+\overrightarrow{BC}=\overrightarrow{AC}$
\vspace{2mm}

\noindent By  axiom A3.2 (invariance under $\sim$ of addition of arrows) the definition is correct because it does not depend on the choice of arrows which represent vectors.

Vectors will be denoted by letters with ''arrows'', for example, $\overrightarrow{a}$.

\begin{theorem}[\textbf{W3}] Addition of vectors makes  set of all vectors into a commutative group:
	
	\begin{enumerate}
		\item $(\overrightarrow{a}+\overrightarrow{b})+\overrightarrow{c}=\overrightarrow{a}+(\overrightarrow{b}+\overrightarrow{c})$
		\item $\overrightarrow{a}+\overrightarrow{0}=\overrightarrow{a}$
		\item $\overrightarrow{a} + (-\overrightarrow{a})=\overrightarrow{0}$
		\item $\overrightarrow{a}+\overrightarrow{b}=\overrightarrow{b}+\overrightarrow{a}$
	\end{enumerate}
\end{theorem}

\begin{proof} We will represent vectors by the corresponding arrows. Let $\overrightarrow{a}=AB$, $\overrightarrow{b}=BC$ and $\overrightarrow{c}=CD$. By the definitions of addition of vectors, null vector and inverse vector, the claims we want to prove  follow from the corresponding claims about arrows.
	
	\begin{enumerate}
		\item $(AB+BC)+CD=AC+CD=CD$,
		
		$AB+(BC+CD)=AB+BD=CD$
		\item $AB+BB=AB$
		\item $AB+(-AB)=AB+BA=AA$
		\item Commutativity of vector addition is a consequence of the parallelogram law (T4). Let's choose a point  $D$ such that $AD\sim BC$. By the parallelogram law, then  $AB\sim DC$ (Fig.\ref{fig:3}).

		\begin{figure}[h]
			\begin{center}
				\includegraphics[height=3cm]{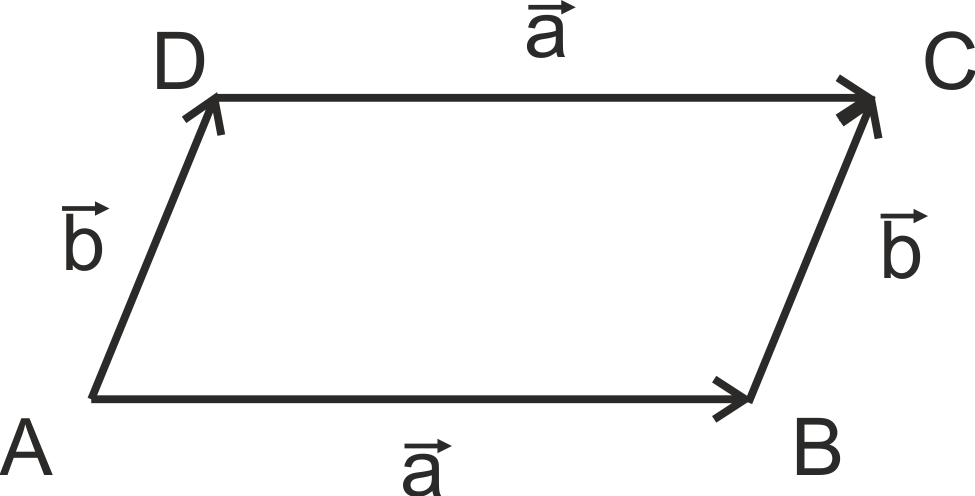}
				\caption{ }
				\label{fig:3}
			\end{center}
		\end{figure}
		
		\noindent Now, we will calculate  the left and right side of the identity we want to prove by "going" along the sides of the parallelogram from the vertex $A$ to the opposite vertex  $C$ in two different ways, via the point $B$ and via  the point $D$:
		
		$\overrightarrow{a}+\overrightarrow{b}=\overrightarrow{AB}+\overrightarrow{BC}=\overrightarrow{AC}$
		
		$\overrightarrow{b}+\overrightarrow{a}=\overrightarrow{AD}+\overrightarrow{DC}=\overrightarrow{AC}$
		
		\noindent Hence, $\overrightarrow{a}+\overrightarrow{b}=\overrightarrow{b}+\overrightarrow{a}$.
	\end{enumerate}\end{proof}
	
	Since the  multiplication of an arrow by a number is invariant under  relation  $\sim$ (axiom A5), we can transfer the operation into the \textbf{multiplication of a vector by a number}:
	
	\vspace{2mm}
	$\lambda\cdot \overrightarrow{AB}=\overrightarrow{\lambda\cdot AB}$
	\vspace{2mm}
	
	Because of the invariance, we can transfer all the properties of the multiplication of an arrow by a number into the properties of the multiplication of a vector by a number:
	
	\begin{theorem}[\textbf{W4}] Multiplication of a vector by a number has the following properties:
		
		\begin{enumerate}
			\item $1\cdot\overrightarrow{a}=\overrightarrow{a}$
			\item $(\lambda +\mu)\cdot\overrightarrow{a}=\lambda\cdot\overrightarrow{a}+\mu\cdot\overrightarrow{a}$
			\item $\lambda\cdot(\mu\cdot \overrightarrow{a})=(\lambda\cdot\mu)\cdot \overrightarrow{a}$
			\item $\lambda\cdot(\overrightarrow{a}+\overrightarrow{b})=\lambda\cdot\overrightarrow{a}+\lambda\cdot\overrightarrow{b}$
		\end{enumerate}
	\end{theorem}
	
	\begin{proof} Let $\overrightarrow{a}=AB$ and $\overrightarrow{b}=BC$
		
		\begin{enumerate}
			\item $1\cdot\overrightarrow{a}=\overrightarrow{1\cdot AB}=\overrightarrow{AB}=\overrightarrow{a}$
			\item $(\lambda +\mu)\cdot\overrightarrow{a}=(\lambda +\mu)\cdot\overrightarrow{AB}=\overrightarrow{(\lambda +\mu)\cdot AB} = \overrightarrow{\lambda\cdot AB +\mu\cdot AB}=\overrightarrow{\lambda\cdot AB}+\overrightarrow{\mu\cdot AB}=\lambda\cdot\overrightarrow{AB}+\mu\cdot \overrightarrow{AB}=\lambda\cdot\overrightarrow{a}+\mu\cdot \overrightarrow{a}$
			\item $\lambda\cdot(\mu\cdot \overrightarrow{a})=\lambda\cdot(\mu\cdot \overrightarrow{AB})=\lambda\cdot\overrightarrow{\mu\cdot AB}=\overrightarrow{\lambda\cdot(\mu\cdot AB)}=\overrightarrow{(\lambda\cdot\mu)\cdot AB}=(\lambda\cdot\mu)\cdot \overrightarrow{AB}=(\lambda\cdot\mu)\cdot \overrightarrow{a}$
			\item $\lambda\cdot(\overrightarrow{a}+\overrightarrow{b})=\lambda\cdot(\overrightarrow{AB}+\overrightarrow{BC})=
			\lambda\cdot\overrightarrow{AB+BC}=\overrightarrow{\lambda\cdot(AB+BC)}=\overrightarrow{\lambda\cdot AB+\lambda\cdot BC}=
			\overrightarrow{\lambda\cdot AB}+\overrightarrow{\lambda\cdot BC}=\lambda\cdot \overrightarrow{AB}+\lambda\cdot \overrightarrow{BC}=\lambda\cdot \overrightarrow{a}+\lambda\cdot \overrightarrow{b}$
		\end{enumerate}
	\end{proof}
	
	\label{Wa}All previously established  propositions with label W say that the space of vectors   $\overrightarrow{S}$ together with operations of addition of vectors and multiplication of a vector by a number has the structure of vector space, and that  the space of points   $S$ together with the space of vectors $\overrightarrow{S}$ and  the mapping $AB \mapsto \overrightarrow{AB}$, that is to say, the structure
	
	\vspace{2mm}
	$(S,\overrightarrow{S},\ ^{\overrightarrow{ \hspace{4mm} }},+,\cdot )$, \ where \ $\ ^{\overrightarrow{ \hspace{4mm} }}: S^2 \rightarrow \overrightarrow{S},\  +:\overrightarrow{S}^2 \rightarrow \overrightarrow{S}, \  \cdot :\mathbb{R}\times \overrightarrow{S} \rightarrow \overrightarrow{S}$,
	\vspace{2mm}
	
	\noindent is Weyl's structure of affine space.  W propositions  are precisely  Weyl's axioms for the affine  part of Euclidean geometry. Conversely, starting  from Weyl's structure of affine space, for which  W propositions are valid, we could define in a standard way the equivalence of arrows, addition of arrows and multiplication of an arrow by a number, and  prove that for such a defined structure
	
	\vspace{2mm}
	$(S,\sim ,+,\cdot)$, \  where \ $S\neq\emptyset, \ \sim \subseteq S^2\times S^2,\ + :S^2\times S^2 \rightarrow S^2  ,\  \cdot : \mathbb{R}\times S^2 \rightarrow S^2$
	\vspace{2mm}

	\noindent all  $A$ propositions are valid. Hence, 
	
	\vspace{2mm}
	\noindent \textit{ A axioms  are equivalent to  W axioms in this affine layer of  Euclidean geometry}.
	\vspace{2mm}
	
	Let's note that in the structure $(S,\sim ,+,\cdot)$ addition $+$ is defined and multiplication $\cdot$ can be defined if we choose to introduce  multiplication of a vector gradually, first by natural numbers,  then by integers and rational numbers, and finally by real numbers. It means that there is essentially only one primitive term, the relation of equivalence $\sim$  between ordered pairs of points.
	
	In the next section we will use a few propositions of affine geometry concerning parallelograms and  projections. It is well known that the propositions follow from W (and so from A) axioms and they will not be proved here.
	
\begin{figure}[h]
		\begin{center}
			\includegraphics[height=4cm]{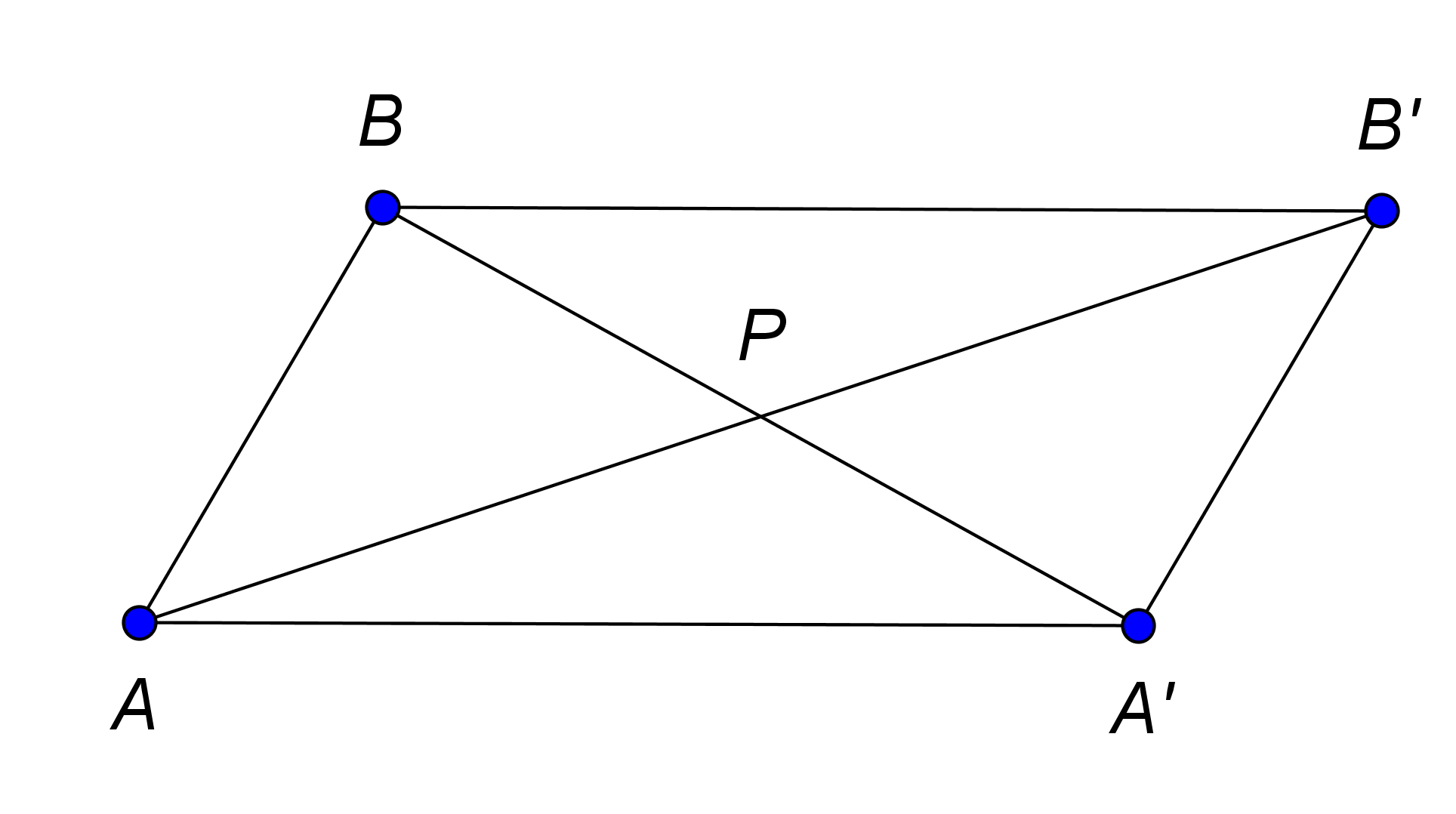}
			\caption{ }
			\label{fig:6}
		\end{center}
	\end{figure}

	We will define a parallelogram in the usual way,  as a quadrilateral $ABCD$ such that $AB\sim DC$. We will need the following result about parallelograms:
	
	\vspace{2mm}
	\noindent \textbf{characterization of a parallelogram by diagonals}. \textit{A quadrilateral is a parallelogram if and only its diagonals bisect each other} (Fig.\ref{fig:6}).
	
	\vspace{2mm}
	$AB\sim A'B' \ \leftrightarrow\ AP\sim PB' \ \land\ BP\sim PA'$
	\vspace{2mm}
	

\begin{figure}[h]
		\begin{center}
			\includegraphics[height=4cm]{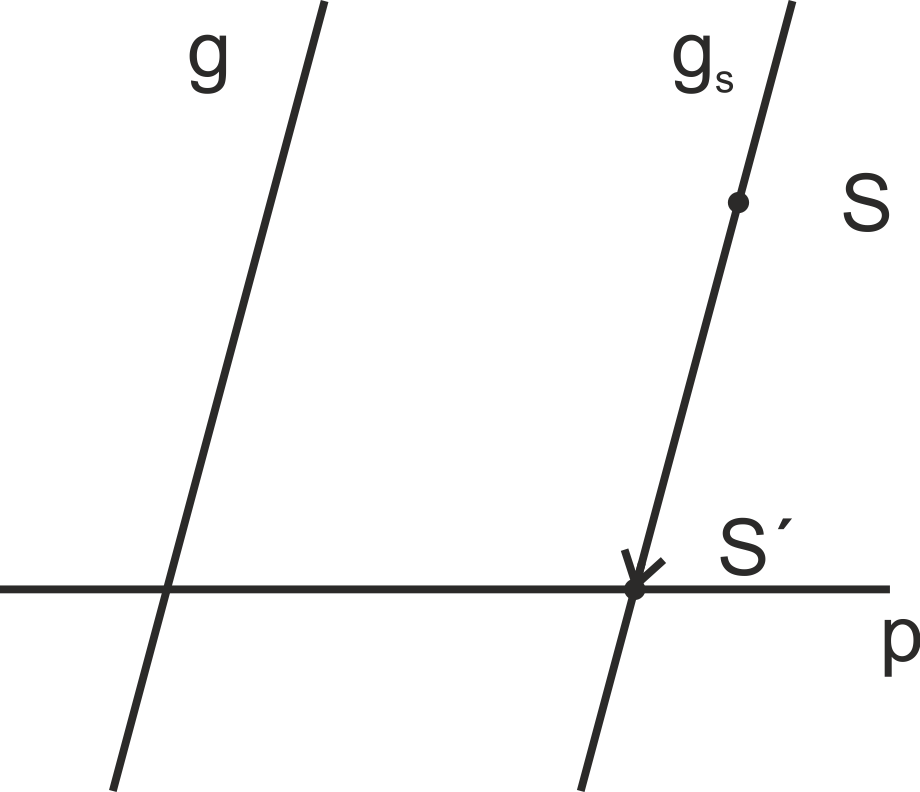}
			\caption{ }
			\label{fig:6A}
		\end{center}
	\end{figure}

	\textbf{The definition of a projection} is based on the following proposition. \textit{In a fixed plane, for every (straight) line  $g$  and point $S$ there is a unique line  $g_S$ through $S$ which does not intersect  $g$ or is equal to $g$ (when $S$ lies on $g$). The line $g_S$ intersects every line $p$ which intersects  $g$ in a unique point  $S'$. The mapping $S \mapsto S'$  maps every point of the plane onto a point of the line $p$. This mapping is termed \textbf{projection}  $P_{g,p}$ onto the line $p$ generated by the line  $g$ }(Fig.\ref{fig:6A}).

	We will need  two \textbf{properties of projection}:
	
	\begin{enumerate}
		
		\item \textit{Projection maps the sum of arrows onto the sum of the projections of the arrows} (Fig.\ref{fig:4}):

		\begin{figure}[h]
			\begin{center}
				\includegraphics[height=4cm]{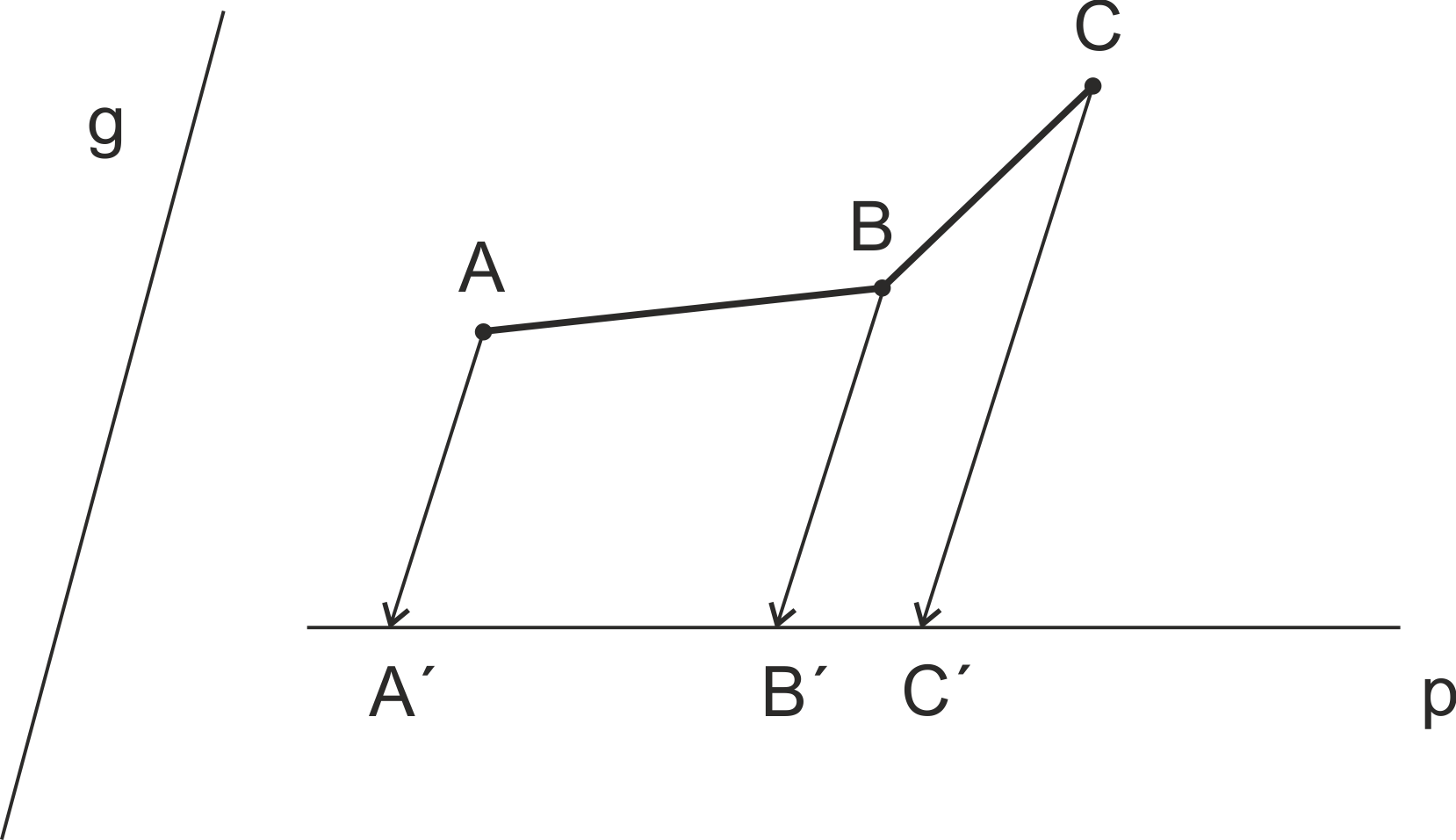}
				
				$P_{g,p}(AB+BC)=P_{g,p}(AB)+P_{g,p}(BC)$
				\caption{ }
				\label{fig:4}
			\end{center}
		\end{figure}	
%
		
		\item \textit{Projection maps  a $\lambda$ times longer arrow into a  $\lambda$ times longer projection of the arrow} (Fig.\ref{fig:5}):
		
	\begin{figure}[h]
		\begin{center}
			\includegraphics[height=4cm]{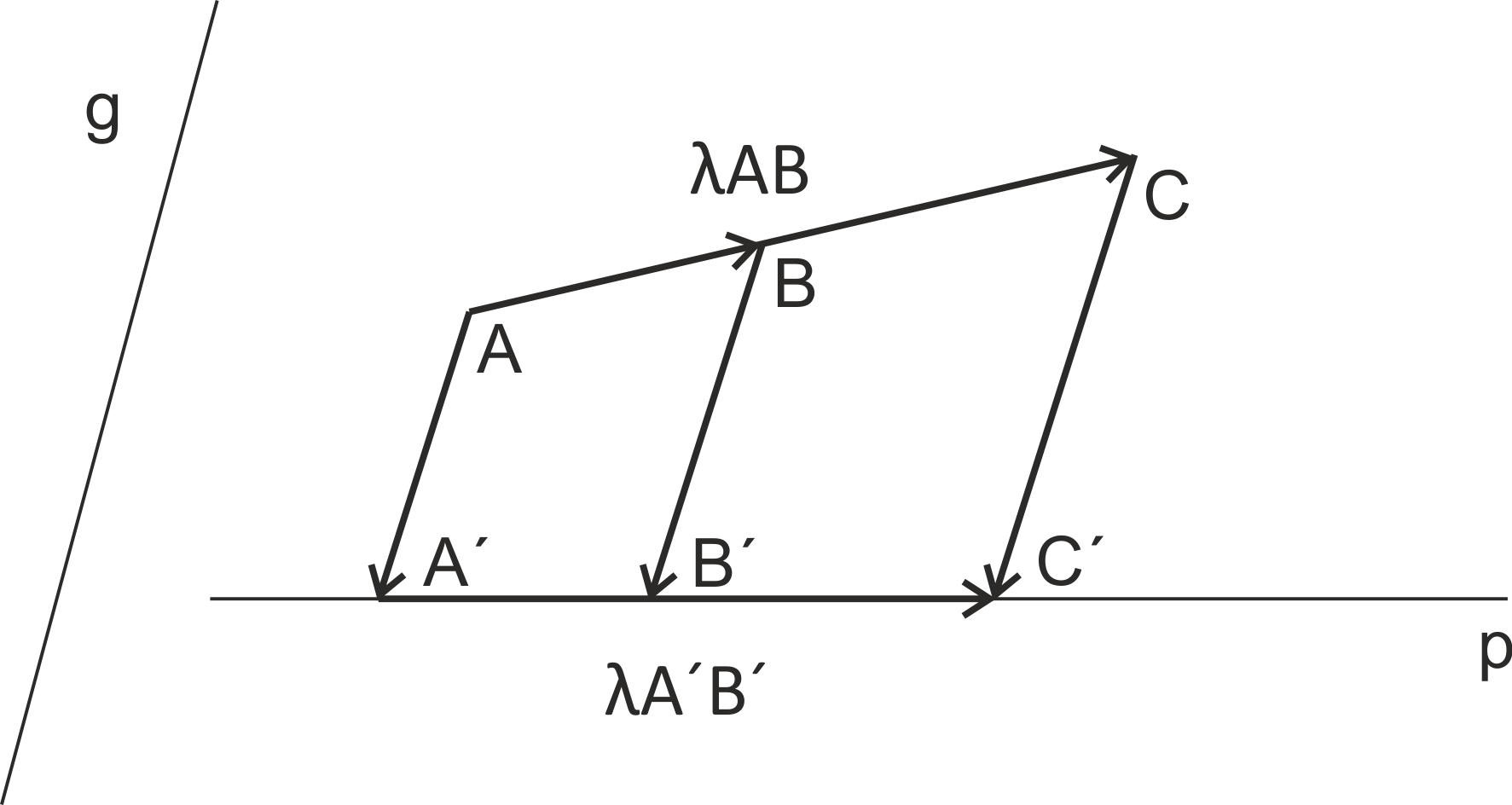}
			
			$P_{g,p}(\lambda AB)=\lambda P_{g,p}(AB)$
			
			\caption{ }
			\label{fig:5}
		\end{center}
	\end{figure}
%
		
	\end{enumerate}

	\section{Length}
	
	A basic geometric measure is the measure of distance between  points $A$ and $B$, the function $|\hspace{3mm}|:S^2 \rightarrow \mathbb{R}$. This is the next and the final primitive term of the theory I develop here. A real number  $|AB|$  will be termed  the \textbf{length} of the arrow $AB$  or  \textbf{distance} from the point $A$ to the point  $B$.

	\noindent By the homogeneity of space the length of an arrow must be invariant under equivalence relation $\sim$:
	
	\begin{axiom}[\textbf{A8}]  $AB\sim CD \ \rightarrow\ |AB|=|CD|$. \end{axiom}
	
	\noindent By the very idea of measuring distance:
	
	\begin{axiom}[\textbf{A9.1}]  $|AA|=0$. \end{axiom}
	
	Every point $B\neq A$ determines a direction in which we can go from $A$. Because of the isotropy of space, the algebraic sign of  distance must be always the same -- distance must be always negative or always positive or always zero. The zero case gives a trivial measure which does not make any difference between arrows, so, it is a useless measure. Thus, the two other possibilities remain. Technically speaking they are mutually equivalent choices, but by the very idea of measuring it is natural to choose a positive algebraic sign:
	
	\begin{axiom}[\textbf{A9.2}]   $B\neq A  \ \rightarrow\ |AB|>0$. (positive definiteness) \end{axiom}
	
	By the isotropy of space we also have:
	
	\begin{axiom}[\textbf{A9.3}]   $|AB|=|BA|$. \end{axiom}
	
	For every direction from a point  $A$ determined with a point $B\neq A$ we already have a measure of distance. If we take $AB$ as a unit of measure, than we can take the number $\lambda > 0$ as a measure of distance of $AC$ where $AC=\lambda AB$. Note that such a choice of measure along every direction need not  be isotropic. However, along every direction the measure of distance $A,B \mapsto |AB|$ must be in accordance with this $\lambda$ measuring (although it must be more than this):
	
	\begin{axiom}[\textbf{10}]    $|\lambda AB|=\lambda |AB|$, for $\lambda >0$, \end{axiom}
	
	We can express axioms  A9.1, A9.3 and A10 in a uniform way by the next equivalent proposition:
	
	\begin{theorem}[\textbf{7}] (compatibility of distance with multiplication) 
		
		$|\lambda AB|=|\lambda ||AB|$, for every real number $\lambda$. \end{theorem}
	
	\begin{proof} By A10, the claim is valid for  $\lambda >0$. By A9.1 the claim is valid for $\lambda =0$. By A9.3 the claim is valid for $\lambda =-1$. Thus, we must prove the claim for the remaining negative values of  $\lambda$. Let $\lambda <0$. Then $\lambda=-\mu$, where $\mu>0$. Now we calculate using what  we have already proven:
		
		$|\lambda AB|=|(-1)\mu AB|=|-1||\mu AB|=|-1||\mu||AB|=|(-1)\mu||AB|=|\lambda ||AB|$\end{proof}
	
	The description of distance function we have achieved until now enables us to compare distances in a given direction with distances in the opposite direction and with distances in parallel directions. What remains is to solve the main problem: how to compare distances along arbitrary directions in an isotropic way.  Let's take, in a given plane, along every direction from a point $S$, a point at a fixed distance $r>0$ from $S$. The set of such  points is the \textbf{circle} with center $S$ and radius $r$, $C(S,r)=\{ T: |ST| = r\}$. Let's choose two points $A$ and $B$  on the circle and consider the unique line   $p(A,B)$ through these points (Fig.\ref{fig:7}):
	
\begin{figure}[h]
	\begin{center}
		\includegraphics[height=3.5cm]{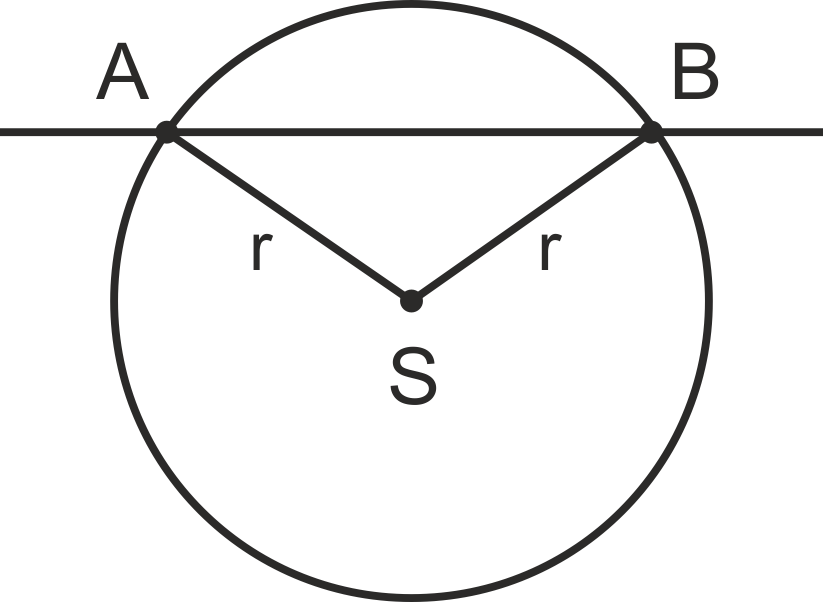}
		\caption{ }
		\label{fig:7}
	\end{center}
\end{figure}
	
	\noindent Let's take an arbitrary point $T$ on the line $p(A,B)$ and consider  how   the distance  $d(T)$ from $T$ to the center $S$ of the circle varies with the choice of $T$.  Thereby, we will use the idea of continuity of space and  of continuity of  function $d(T)$. Because of the isotropy of space, the function $d(T)$ must be symmetric  in the position of $T$  relative to  the points  $A$ and $B$ (directions $SA$ and $SB$). For example, values of the function in the points  $A$ and $B$ are the same (equal to $r$). Also, the function must have the same value in a point we reach when we move  a certain distance from  $A$ to $B$ as well as in a point we reach when we move  the same distance from $B$ to $A$: $d(A+ \lambda \overrightarrow{AB})=d(B+\lambda \overrightarrow{BA})$ (Fig.\ref{fig:8}):
	
	\begin{figure}[h]
		\begin{center}
			\includegraphics[height=4cm]{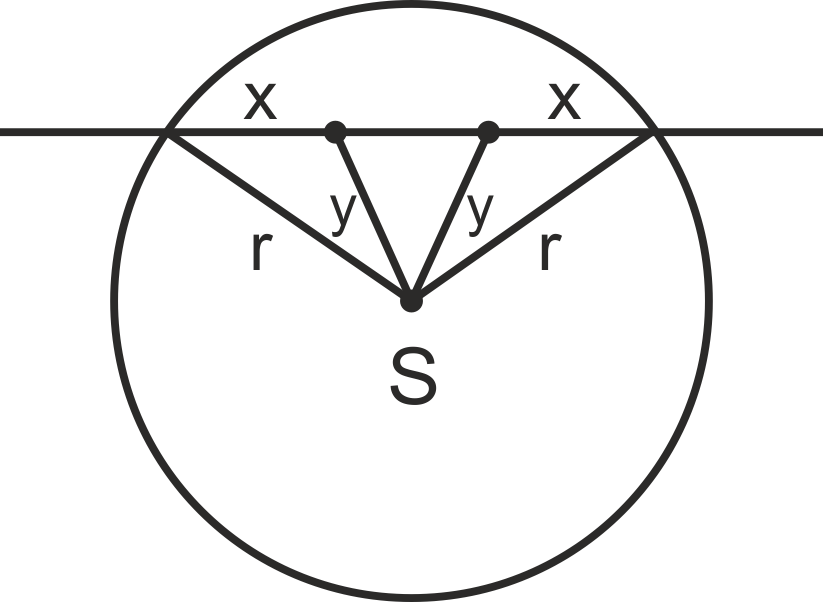}
			\caption{ }
			\label{fig:8}
		\end{center}
	\end{figure}
%
	
	\noindent Because of this symmetry, the function  $d(T)$ must have a local extreme value in the midpoint of $AB$. To determine more precisely the character of the extreme point we will exploit   knowledge of a special case,  when the points   $A$ and $B$  are diametrically opposite on the circle, that is to say, when the center $S$ of the circle lies on  $p(A,B)$. In that case, if we ''move'' a point $T$ from $A$ to $B$  (or from  $B$ to $A$), the distance $d(T)$ from the center $S$ of the circle decreases and it is  smallest in the midpoint ($S$). Furthermore,  if we move $T$ from $A$ in the direction opposite to the direction to $B$  (or from $B$ in the direction opposite to the direction to $A$), the distance increases. Therefore,  the midpoint $S$ is a unique point of the global minimum of the function $d(T)$. If we drag the point $B$ slightly  along the circle into  the point  $B'$, the center  $S$ of the circle will no longer be on the line $p(AB')$, but, because of  continuity, the behaviour  of the function $d(T)$ will remain the same. That is to say, the midpoint  $P$ of $AB'$ will remain a unique  global minimum of the function on the line (Fig.\ref{fig:9}):
	
	\begin{figure}[h]
		\begin{center}
			\includegraphics[height=3.5cm]{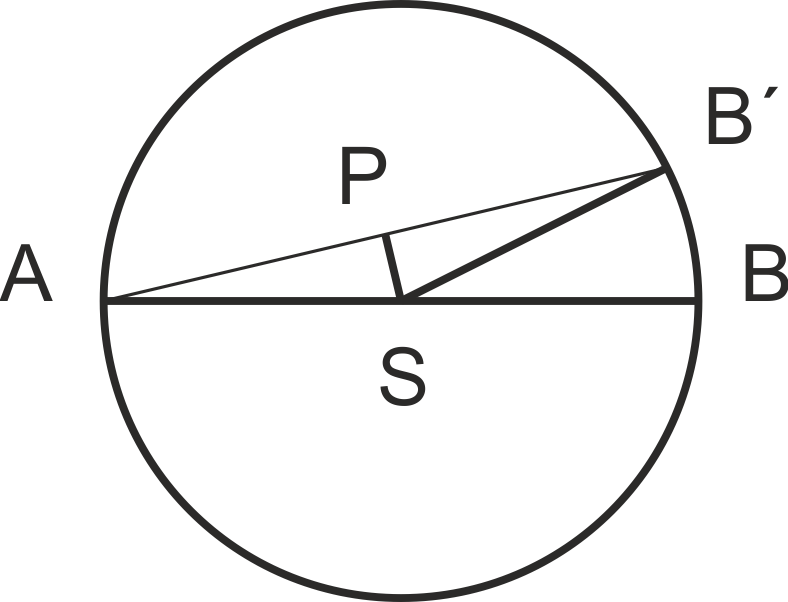}
			\caption{ }
			\label{fig:9}
		\end{center}
	\end{figure}

\begin{figure}[h]
			\begin{center}
				\includegraphics[height=3.5cm]{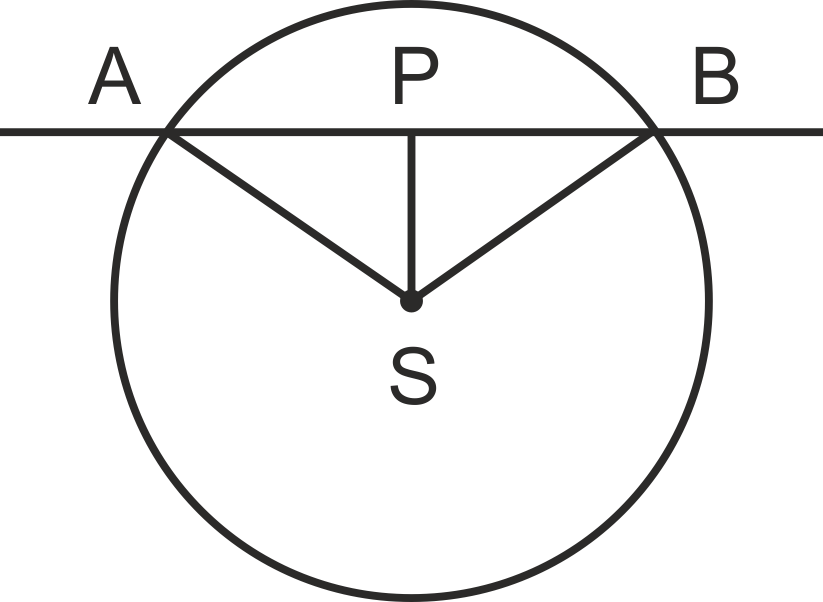}
				\caption{ }
				\label{fig:10}
			\end{center}
		\end{figure}

	\noindent Because of continuity,  for every two points  $A$ and $B'$ on the circle the function $d(T)$ will have a unique global minimum on  line $p(AB')$ exactly in the midpoint of $AB'$:
	
	\begin{axiom}[\textbf{11}]    If a line has two   common points with a circle, points $A$ and $B$, then the midpoint $P$ of $AB$  is the point on the line nearest to the center of the circle. \normalfont{(Fig.\ref{fig:10})} \end{axiom}
	
	From the axiom it follows immediately that a line can not have more than two common points with a circle. If, besides $A$ and $B$, there was a third common point $X$, by A11, the midpoint of $AB$ and the midpoint of  $AX$ would be  points on a line nearest to the center of the circle. Then, by uniqueness of the nearest point,  $AB$ and   $AX$ would have the same midpoint. Thus, we would get a contradiction, that $X=B$.
	
	Let a line  $p$ have exactly one common point with a circle, a point  $A$. If we drag the point $A$  slightly  along the circle  in one direction onto a point $Al$, and in another direction onto a point  $Ad$, then the line $p$ is dragged  onto the line  $p(Al,Ad)$. By  axiom A11 the midpoint  $P$ of $AB$ is the point on $p(Al,Ad)$ nearest to the center of the circle. By continuity of space, the point  $A$ must be the point on  $p$ nearest to the center of the circle (Fig.\ref{fig:11}):
	
\begin{figure}[h]
	\begin{center}
		\includegraphics[height=4cm]{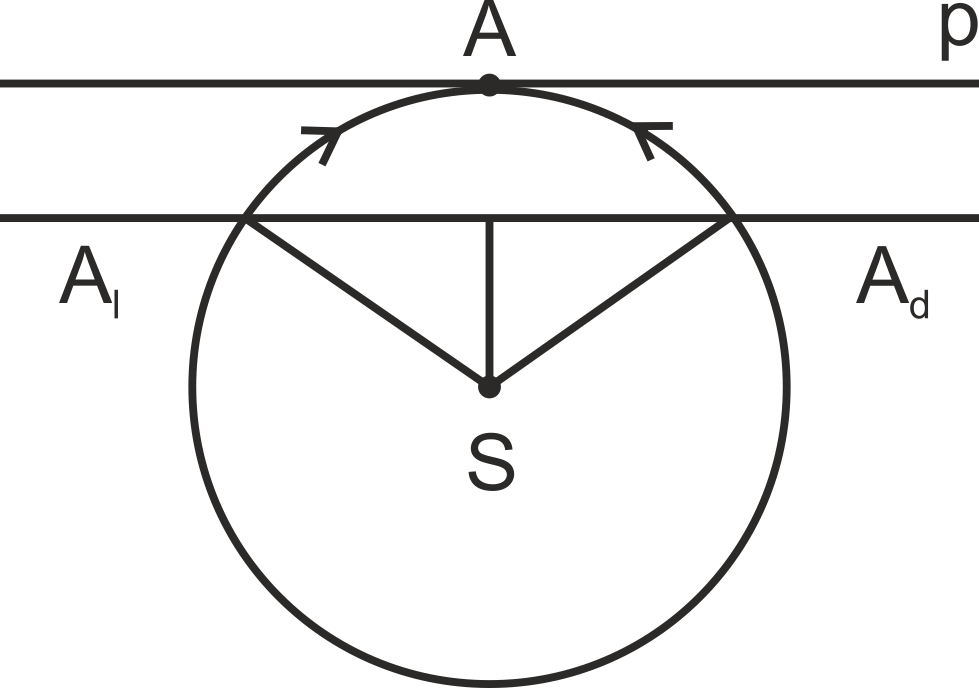}
		\caption{ }
		\label{fig:11}
	\end{center}
\end{figure}
	
\begin{figure}[h]
		\begin{center}
			\includegraphics[height=3.5cm]{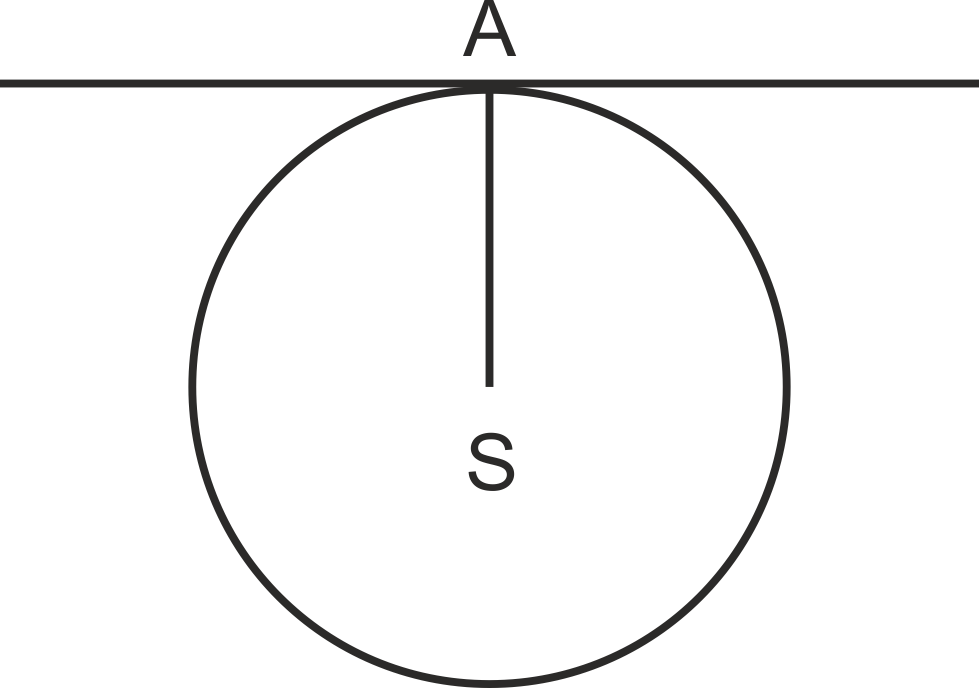}
			\caption{ }
			\label{fig:12}
		\end{center}
	\end{figure}

	\begin{axiom}[\textbf{A12}]   If a line has exactly one common point with a circle, then the common point is the point on the line nearest to the center of the circle. \normalfont{(Fig.\ref{fig:12})}\end{axiom}

	\begin{theorem}[\textbf{7}] For every point  $S$ not on a line  $p$  there is a unique point  $P$ on $p$ which is the point on $p$ nearest to $S$. \end{theorem}
	
	\begin{proof} Let $A$ be a point on a line $p$ and let its distance to $S$ be $r$. Thus, $A$ is a common point of the line $p$ and the circle with center $S$ and radius $r$. If the line $p$ does not have another  common point with the circle, then, by A12,  $A$ is a point on $p$ nearest to  $S$.  If line $p$ has another  common point $B$ with the circle, then, by A11, the midpoint $P$ of $AB$ is the point on $p$ nearest to  $S$.  Thus, we have proved the existence of the nearest point. The  nearest point is unique by its very definition. \end{proof}
	
	The direction from a point not on a line to its nearest point on the line is a perpendicular direction to the line. Formally, we define that a line  $a$ is \textbf{perpendicular} to a line $b$, in symbols $a\bot b$, if $a$ intersects  $b$ and there is a point  $S$ on $a$ which is not on $b$ such that the intersection of $a$ and $b$ is the point on $b$ nearest to $S$.
	
	\begin{theorem}[\textbf{8}] For lines $a$ and $b$, if $a\bot b$ then $b\bot a$.\end{theorem}
	
	\begin{proof} By definition of perpendicularity, there is a point $S$ on $a$ such that  the intersection of $a$ and $b$, the point $O$, is the point on $b$ nearest to $S$.  Let $P$ be a point on  $b$ different from $O$. Since its distance to $S$ is greater than the distance of  $O$ to $S$, there is another point $P'$ on the line $b$ equally distant as $P$ to $S$ (Fig.\ref{fig:13}). On the contrary, by A12, $P$ would be the point on   $b$ nearest to $S$, hence, it would be  $P=O$, which is  in a contradiction of the choice of $P$. Since $O$ is the point on $b$ nearest to $S$, it is, by A11, the midpoint of $PP'$.
		
	\begin{figure}[h]
		\begin{center}
			\includegraphics[height=3.5cm]{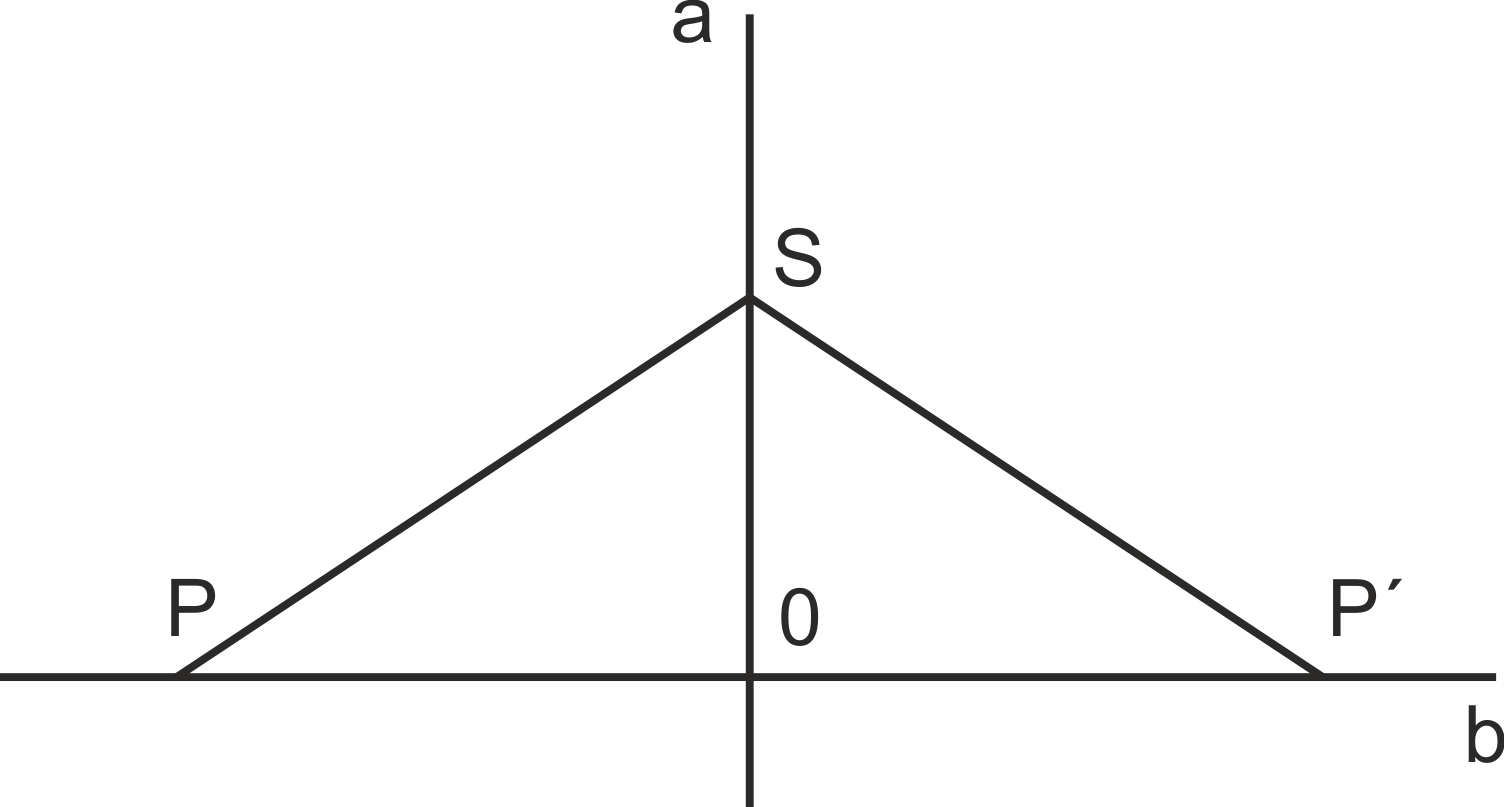}
			\caption{ }
			\label{fig:13}
		\end{center}
	\end{figure}
	

		\noindent  By the uniqueness of the translation of arrows law T3. there is  a point $S'$ on the line $a$  such that $SO\sim OS'$ (Fig.\ref{fig:14}). Thus, the point $O$ is the midpoint of $SS'$, too. Therefore, diagonals of the quadrilateral $SPS'P'$ bisect each other.
		
	\begin{figure}[h]
		\begin{center}
			\includegraphics[height=4cm]{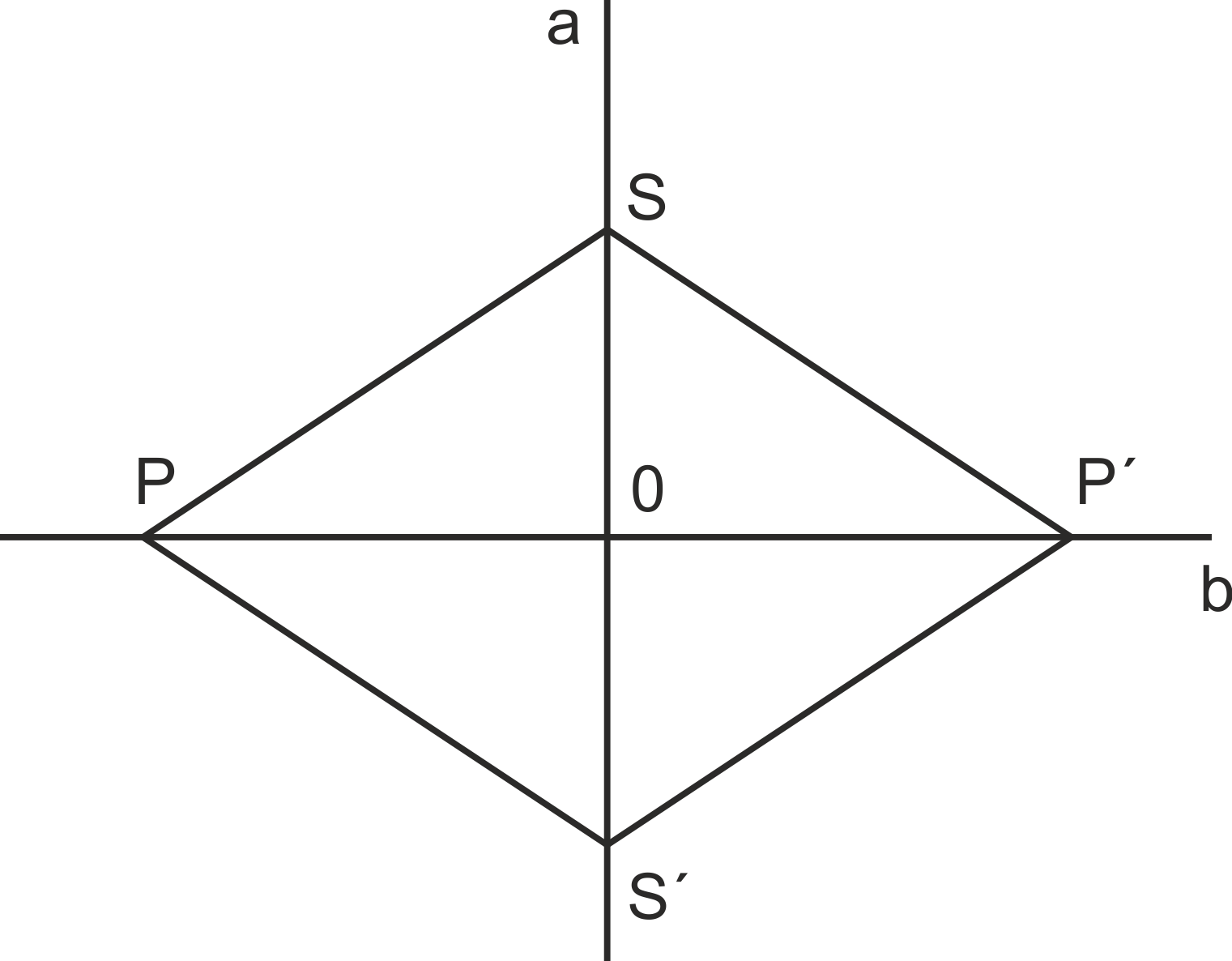}
			\caption{ }
			\label{fig:14}
		\end{center}
	\end{figure}

		\noindent    By  the theorem of affine geometry on the characterization of a parallelogram by diagonals, the quadrilateral is a parallelogram. Since it is a parallelogram, the opposite sides have the same length. However, by the choice of the point $P'$, the neighbouring sides have the same length, too. Therefore, all sides of the quadrilateral have the same length. Thus,  $PS$ and $PS'$ have the same length. Hence, by A11, the point $O$ is the point on $a$ nearest to $P$. By the definition of perpendicularity, it means that the line $b$ is perpendicular to the line $a$.\end{proof}

	\begin{theorem}[\textbf{9}]  A line  $b$ is perpendicular to a line $a$ if and only if they have a common point $O$ and for every point $P$ on line $b$ the intersection  $O$ is the point on  the line  $a$ nearest to $P$. \end{theorem}
	
	\begin{proof} One direction is trivial: if for every point   $P$ on line $b$ the intersection $O$ is the  point on line $a$ nearest to $P$ then that is true for a particular $P$ on $b$, so, by definition of perpendicularity,  $b\bot a$. The opposite direction is the immediate consequence of the previous theorem and its proof. Namely, by the theorem, from $b\bot a$  follows $a\bot b$. In the proof of the previous theorem, from that assumption it is proved  that for every point $P$ on the line $b$ the intersection $O$ is the point on the line  $a$ nearest to $P$. \end{proof}
	
	\begin{corollary}[\textbf{10}] (the triangle inequality theorem):  $|AB|+|BC|\geq  |AC|$. The equality is valid if and only if the point   $B$ is on $AC$.\end{corollary}

	\begin{proof} We will prove only the most important case (Fig.\ref{fig:15}), when the point $B$ is not on line $p(A,C)$ and the line through $B$ perpendicular to $p(A,C)$  intersects $p(A,C)$ at a point $B'$ which is on $AC$ (by theorem 7, the perpendicular line exists).
		
	\begin{figure}[h]
		\begin{center}
			\includegraphics[height=3.5cm]{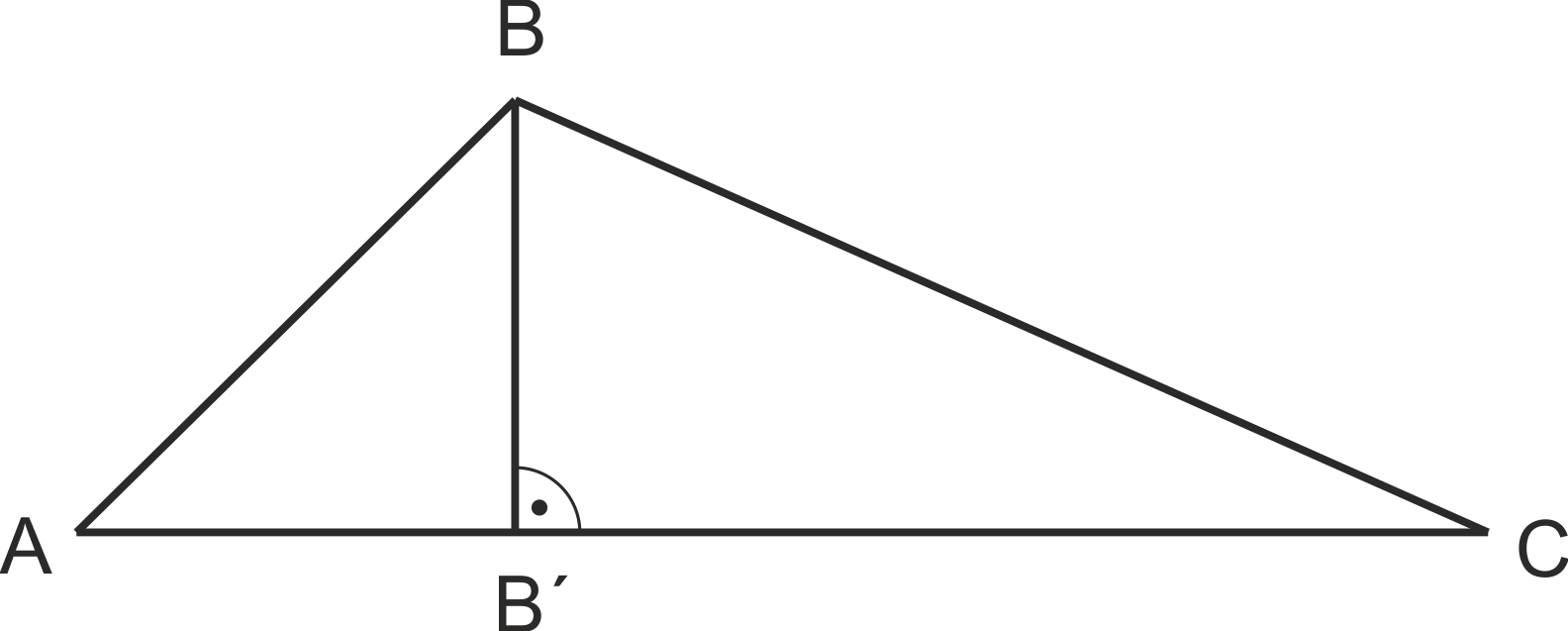}
			\caption{ }
			\label{fig:15}
		\end{center}
	\end{figure}
		
		\noindent By theorem 8, the line $p(A,C)$ is perpendicular to  the line  $p(B,B')$. By theorem 9, $B'$ is the point on $p(B,B')$ nearest to points  $A$ and $C$. Hence, $|AB|>|AB'|$ and $|BC|>|B'C|$. Thus,   $|AB|+|BC|> |AB'|+|B'C| = |AC|$.\end{proof}

	\noindent Theorem 9  enables us to show  that for every point and every line there is a unique line  through the point  perpendicular to the given line .
	
	\begin{theorem}[\textbf{11}] For a point $S$  not on a line $p$ there is a unique line through the point $S$ perpendicular to the line $p$.\end{theorem}
	
	\begin{proof} By theorem 7, there is a point $O$  on the line $p$ nearest to $S$. By the definition of  perpendicularity, the line   $a=p(S,O)$ is perpendicular to $p$. Thus, we have proved the existence of the perpendicular line. Let $S$ lie on another  perpendicular line $b$.   By theorem 10,  the intersection  $O'$ of the line $b$ and the line  $p$ is the point on  $p$ nearest to $S$. By the uniqueness of the nearest point, $O'=O$. Thus, $b$ as well as  $a$ contains $S$ and $O$, hence $b=a$.\end{proof}
	
	\begin{theorem}[\textbf{12}] For every point $P$ on a line $p$ there is a unique line through $P$ perpendicular to $p$.\end{theorem}
	
	\begin{proof} To show the existence of the perpendicular line let's choose a point  $S$ not on  $p$. If the intersection $O$ of  the line through $S$ perpendicular to $p$  with line $p$ is just the point $P$ then we have found the perpendicular line. If it is not so, let $S'$ be a point such that  $SS'\sim OP$ (Fig.\ref{fig:16}):
		
			\begin{figure}[h]
				\begin{center}
					\includegraphics[height=3.5cm]{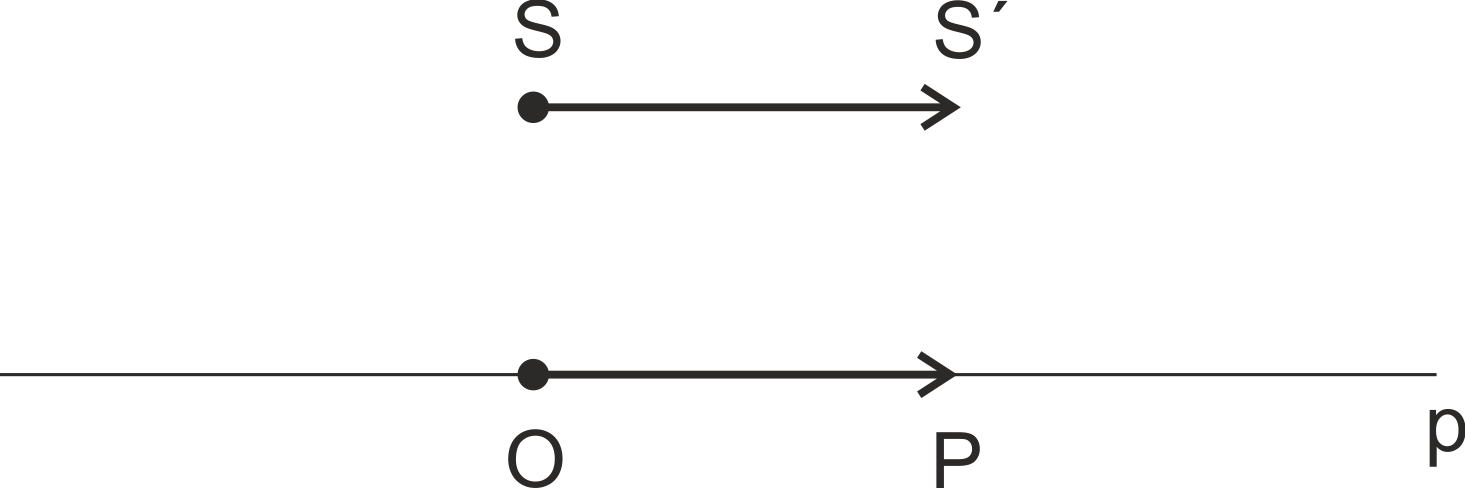}
					\caption{ }
					\label{fig:16}
				\end{center}
			\end{figure}

		\noindent Because of the homogeneity invariance of length (axiom A8), since $O$ is a point on $p$ nearest to $S$, the point  $P$ is the point on $p$ nearest to  $S'$. Hence, $P$ is on  the line through $S'$ perpendicular to line $p$.
		
		\noindent To prove the uniqueness of the perpendicular line, let's suppose the contrary, that there are two lines through the point $P$ perpendicular to the line $p$. Let $S$ be a point on one of  perpendicular lines, the line $a$, but not on $p$ and not on the other perpendicular line, the line $b$. By affine geometry, there is a unique line through $S$ parallel to  $p$. Let $S'$ be the intersection point of the parallel line   with another perpendicular line $b$ (Fig.\ref{fig:17}):
		
	\begin{figure}[h]
		\begin{center}
			\includegraphics[height=4cm]{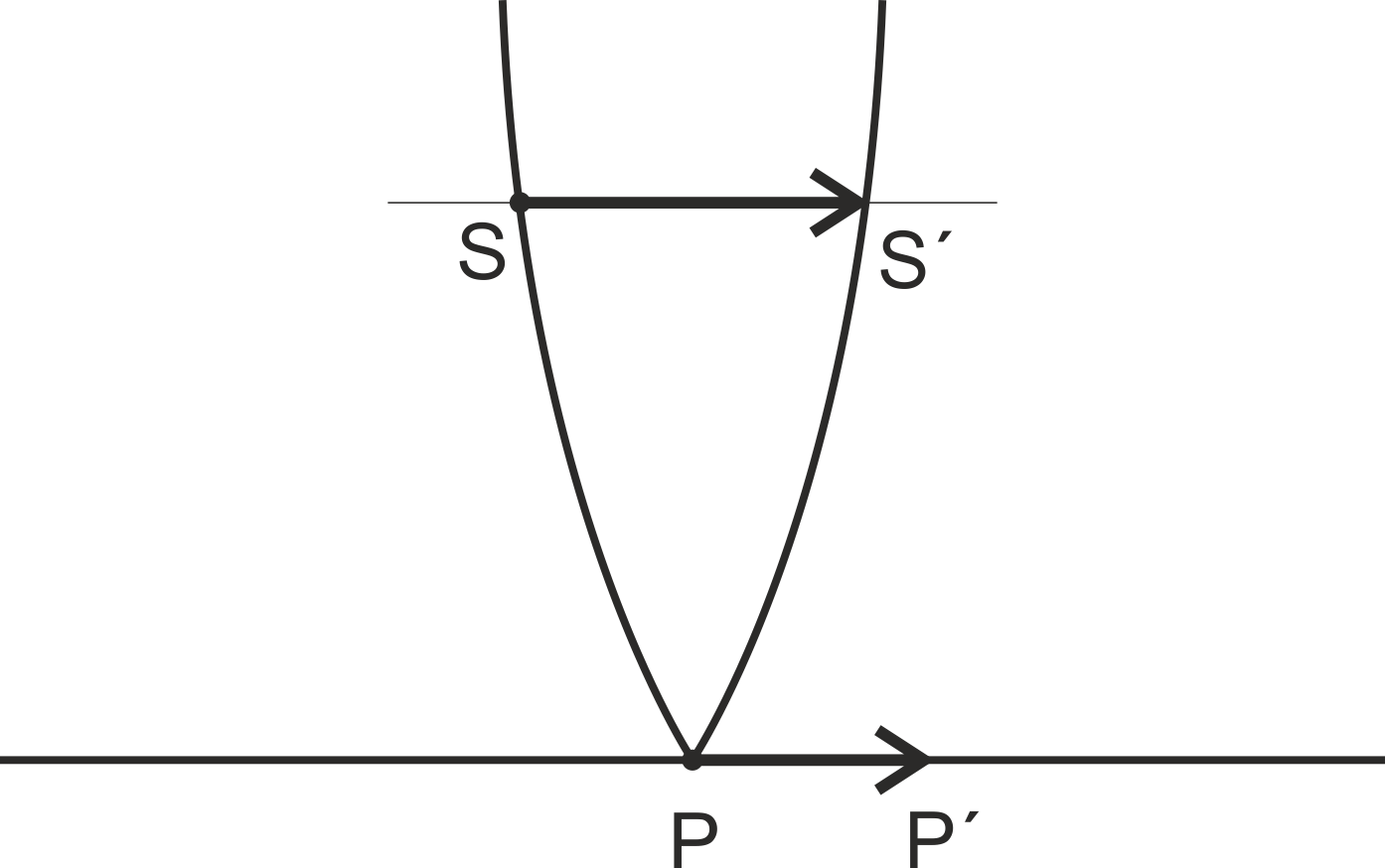}
			\caption{ }
			\label{fig:17}
		\end{center}
	\end{figure}
		
		\noindent Let $P'$ be a point on line $p$ such that $PP'\sim SS'$ Let's observe that $P'\neq P$. Because of the homogeneity invariance of length, since $P$ is the point on $p$ nearest to point $S$, so the  point $P'$ is a point on $p$ nearest to point $S'$. However, since $P$ is on the line through $S'$ perpendicular to $p$, $P$ is the nearest point. Thus, we get a contradiction, that $P'= P$. Therefore, there is no more than one line through the point $P$ perpendicular to the line $p$.\end{proof}

	We will define the \textbf{orthogonal projection} of a point $S$ to a line $p$ to be the intersection of the line $p$ and the unique line through  $S$ perpendicular to  $p$, that is to say, the point on line $p$ nearest to point $S$. The next theorem says that it is indeed a projection.
	
	\begin{theorem}[\textbf{13}] All lines perpendicular to a line  $p$ in a given plane containing $p$ are mutually parallel.\end{theorem}
	
	\begin{proof} In a given plane,  lines $a$ and $b$ are parallel if they have  no common points, or they are the same lines. Let $a$ and $b$ be lines perpendicular to a line $p$ in a given plane containing $p$. If $a$ and $b$ have no common points they are parallel. Let them have a common point  $S$. By theorems  11 (if $S$ is not on $p$) and 12 (if $S$ is on $p$) there is a unique line through $S$  perpendicular to $p$. This means that $a=b$. So, in this case, $a$ and $b$ are parallel, too.\end{proof}
	
	\noindent Therefore, in a given plane, parallel lines generated by a line $g$ perpendicular to a line  $p$ are precisely lines perpendicular to the line $p$.
	
\begin{figure}[h]
		\begin{center}
			\includegraphics[height=3.5cm]{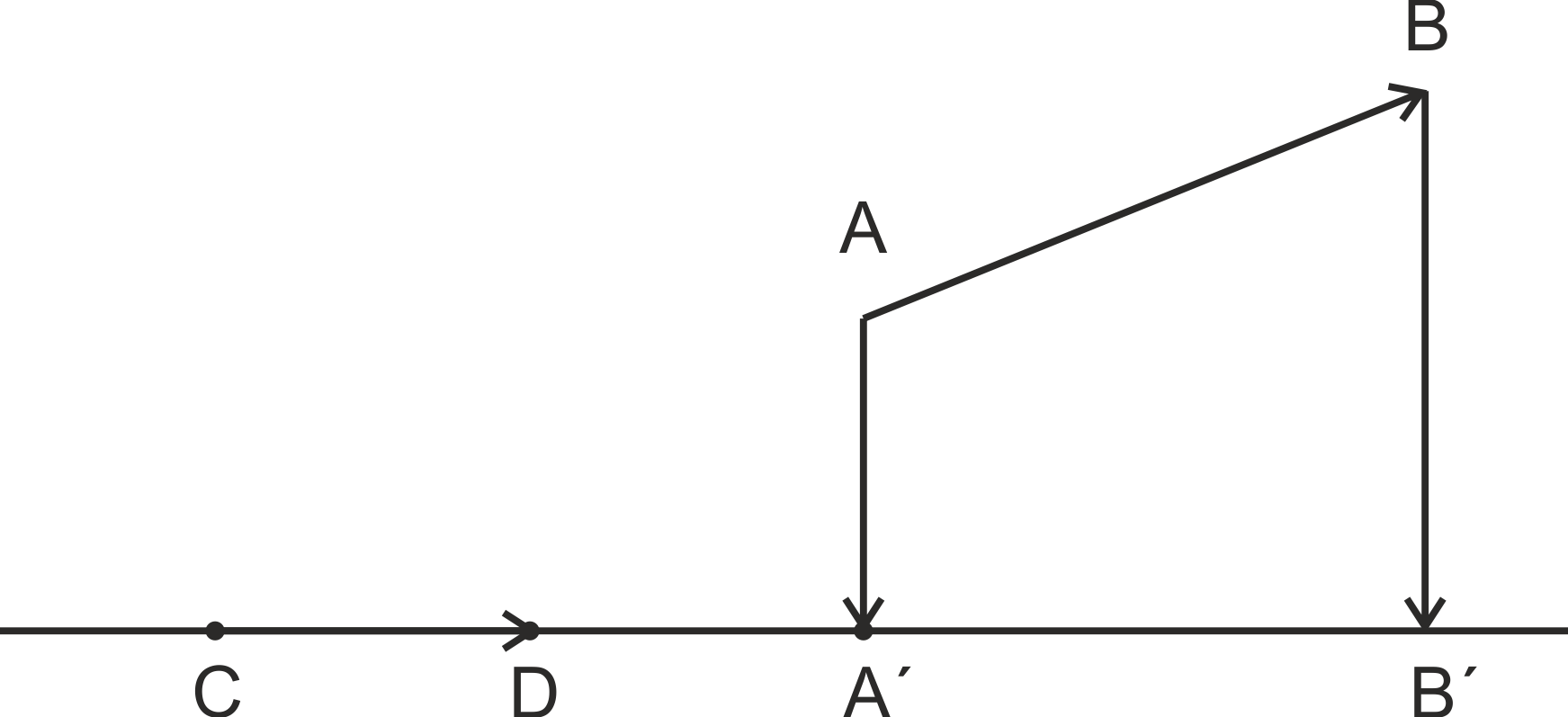}
			\caption{ }
			\label{fig:18}
		\end{center}
	\end{figure}

	Orthogonal projection enables us to define  the scalar orthogonal projection of an arrow onto another arrow.  Let $C\neq D$, and let points $A$ and $B$ be orthogonally projected on line $p(CD)$ into points $A'$ and $B'$ (Fig.\ref{fig:18}). Then  $A'B'\sim \alpha CD$ for some real number $\alpha$.
	
	We define the \textbf{scalar orthogonal projection} of the arrow  $AB$ onto the arrow $CD$ to be the number $\alpha |CD|$. In simpler terms, it is just the $\pm$ length of the orthogonal projection of the arrow $AB$ onto the line $p(CD)$, where the sign is $+$ if the projection is in the direction of $CD$, $-$ otherwise. In the extreme case of null arrow $CC$ it is convenient to take zero for the value of the scalar projection on $CC$. We will denote  $AB_{CD}$ as the scalar projection of  $AB$ onto $CD$.
	
\begin{figure}[h]
	\begin{center}
		\includegraphics[height=4.5cm]{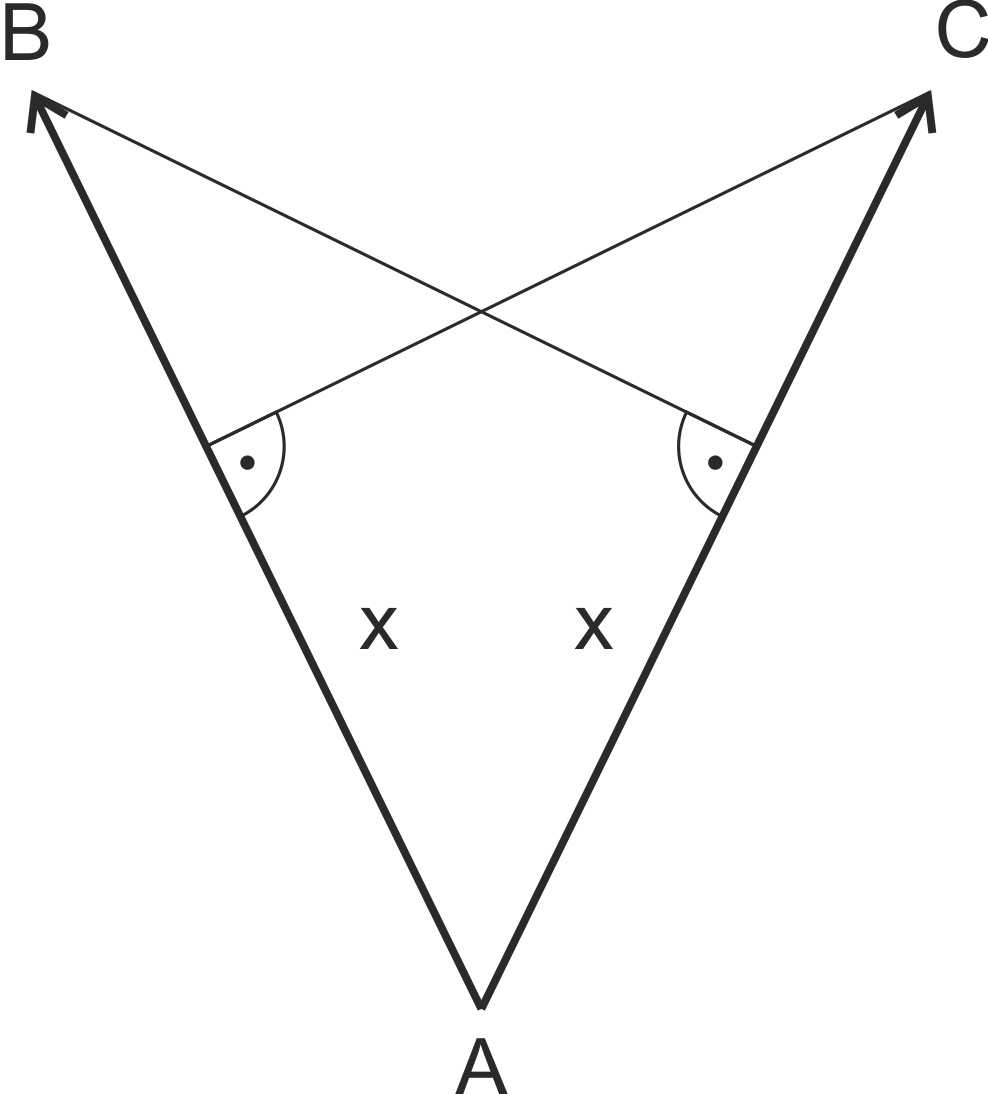}
		\caption{ }
		\label{fig:19}
	\end{center}
\end{figure}

	For two equally long arrows with the same initial point, because of the isotropy of space,  the scalar projection of the first arrow on the second arrow must be the same as the scalar projection of the second arrow on the first arrow. This is the content of the final axiom:
	
	\begin{axiom}[\textbf{A13}]   $|AB|=|AC| \ \rightarrow \ AB_{AC}=AC_{AB}$. \normalfont{(Fig.\ref{fig:19})}\end{axiom}
	
	Now, we can define the \textbf{scalar product of arrows} $AB$ and $CD$. It is the product of the scalar projection of the arrow $AB$ onto $CD$ and  the length of the arrow  $CD$. More formally:
	
	$AB\cdot CD = AB_{CD}\cdot |CD|$
	
\noindent	This operation is invariant under $\sim$ relation:

	\begin{theorem}[\textbf{14}]  $A'B'\sim AB \ \wedge \ C'D'\sim CD \ \rightarrow\ A'B'\cdot C'D'=AB\cdot CD$. \end{theorem}
	
	\begin{proof} It is the consequence of the properties of projection that (i) equivalent arrows are projected onto equivalent arrows and  (ii) projections generated by the same line   onto parallel lines map an arrow onto two mutually equivalent arrows. \end{proof}
	
	\begin{theorem}[\textbf{15}] (properties of the scalar product of arrows)
		
		1. $A\neq B \ \rightarrow AB\cdot AB > 0$
		
		2. $(AB+A'B')\cdot CD = AB\cdot CD + A'B'\cdot CD$
		
		3.  $(\lambda AB)\cdot CD = \lambda (AB\cdot CD)$ \hspace{1cm} $AB \cdot (\lambda CD ) = \lambda (AB\cdot CD)$

		4.  $AB\cdot CD = CD\cdot AB$
		
	\end{theorem}
	
	\begin{proof} \
		
		1. $AB\cdot AB = AB_{AB}\cdot |AB|= |AB| \cdot |AB| > 0$, because, by A9.2, for $A\neq B$ \  $|AB| > 0$.
		
		2.  For the proof it is crucial that, by theorem 13,  orthogonal projection is a special type of  projection. From affine geometry it is known that for every projection it is valid that  $(AB+A'B')_{CD}= AB_{CD}+A'B'_{CD}$. So,
		
		$(AB+A'B')\cdot CD = (AB+A'B')_{CD}\cdot |CD|= (AB_{CD}+A'B'_{CD})\cdot |CD|= AB_{CD}\cdot |CD|+A'B'_{CD}\cdot |CD| = AB\cdot CD + A'B'\cdot CD$
		
		3. Because orthogonal projection is a special type of  projection and for every projection it is valid that   $(\lambda AB)_{CD}= \lambda (AB_{CD})$, we have
		
		$(\lambda AB)\cdot CD = (\lambda AB)_{CD}\cdot|CD| = \lambda (AB_{CD})\cdot|CD|=\lambda (AB\cdot CD)$
		
		\noindent Thus, the first claim is proved. The second claim is an immediate consequence of  definition of the scalar product of arrows:
		
		$AB \cdot (\lambda CD ) = AB_{\lambda CD}\cdot|\lambda CD|=  |\lambda| AB_{\lambda CD}\cdot|CD| =  \lambda AB_{ CD}\cdot|CD|=\lambda (AB\cdot CD)$
		
		4. We will express  the arrow  $AB$ as a stretched arrow of the arrow $C'D'$ which has  length equal to the length of $CD$,  $AB =\lambda C'D'$. Now, we will apply the previously proved proposition, axiom A13, and invariance under $\sim$ of the scalar product of arrows (theorem 14):

		$AB\cdot CD = (\lambda C'D')\cdot CD = \lambda (C'D'\cdot CD)= \lambda (CD\cdot C'D')  =  CD\cdot (\lambda C'D') =CD\cdot AB$
	\end{proof}

	The \textbf {scalar product of vectors} is defined by arrows which represent vectors:
	
	$\overrightarrow{AB}\cdot \overrightarrow{CD}=AB\cdot CD$
	
	\noindent  Since scalar product of arrows is invariant under  $\sim$  (theorem 14)  the definition is correct, that is to say, it doesn't depend on the choice of arrows representing vectors.
	
	\begin{theorem}[\textbf{W5}]   (properties of the scalar product of vectors)
		
		1. $\overrightarrow{a}\neq \overrightarrow{0} \ \rightarrow \overrightarrow{a}\cdot \overrightarrow{a} > 0$
		
		2. $(\overrightarrow{a}+\overrightarrow{b})\cdot \overrightarrow{c} = \overrightarrow{a}\cdot \overrightarrow{c} + \overrightarrow{b}\cdot \overrightarrow{c}$
		
		3.  $(\lambda \overrightarrow{a})\cdot \overrightarrow{b} = \lambda (\overrightarrow{a}\cdot \overrightarrow{b})$

		4.  $\overrightarrow{a}\cdot \overrightarrow{b} = \overrightarrow{b}\cdot \overrightarrow{a}$\end{theorem}
	
	\begin{proof} These claims, by definition of the scalar product of vectors, are reduced to corresponding claims about the scalar product of arrows, which were proved in the previous theorem.\end{proof}

	\label{We}Proposition W5  about the scalar product of vectors, together with  W propositions about the affine structure of space, form Weyl's axiomatics of Euclidean geometry. Conversely, by Weyl's axioms we can define in a standard way the length of an arrow and deduce all A axioms about length. This means that   A axioms of the structure of set of points   $S$ with  relation $\sim$ between pairs of points,  addition of pairs, multiplication of a pair by a number and  distance function of a pair:
	
	\vspace{2mm}
	$(S,\sim ,+,\cdot , |\hspace{3mm}|)$, \  where
	\vspace{2mm}
	$S\neq \emptyset, \  \sim \subseteq S^2\times S^2,\ + :S^2\times S^2 \rightarrow S^2 , \  \cdot : \mathbb{R}\times S^2 \rightarrow S^2, \ |\hspace{3mm}|:S^2 \rightarrow \mathbb{R}$,
	\vspace{2mm}
	
	\noindent are equivalent to  Weyl's axioms of the corresponding structure of the set of points and the set of vectors together with the operation from pairs of points to vectors, addition of vectors, multiplication of a vector by a number and the scalar product of vectors:
	
	\vspace{2mm}
	$(S,\overrightarrow{S},\ ^{\overrightarrow{ \hspace{4mm} }},+,\cdot ,\cdot )$, \ where 
	\vspace{2mm}
	$S\neq\emptyset, \ \ ^{\overrightarrow{ \hspace{4mm} }}: S^2 \rightarrow \overrightarrow{S},\  +,\cdot :\overrightarrow{S}^2 \rightarrow \overrightarrow{S}, \  \cdot :\mathbb{R}\times \overrightarrow{S} \rightarrow \overrightarrow{S}$
\pagebreak	
	\appendix
	
	\section*{Appendix}
	
	The primitive terms of the system of axioms are  (i) equivalence  of pairs of points (arrows): $AB\sim CD$, with the intuitive meaning that the position of the point $B$ relative to the point $A$ is the same as the position of the point $D$ relative to the point $C$,  (ii) multiplication of a pair of points (an arrow) by a real number: $\lambda , A, B \ \mapsto \ \lambda \cdot AB$, with the intuitive meaning of stretching the arrow and of iterative addition of the same arrow, and (iii) distance between points: $A, B \ \mapsto \ |AB|\in \mathbb{R}$.
	
	\textbf{Axiom} (\textbf{A1}) $\sim$ \textit{is an equivalence relation}. \\
\indent	(by the very idea  \textit{to be in the same relative position})
	
	\textbf{Axiom} (\textbf{A2}) $AB\sim AC  \ \rightarrow \  B=C$. \\
\indent (by the very idea of the relative position of points)
	
	\noindent The definitions od inverse arrow and of adding arrows:
	
	\vspace{2mm}
	\noindent \textbf{inverting arrow}:\hspace{5mm} $AB \ \mapsto \  -AB=BA$
	
	\vspace{2mm}
	\noindent \textbf{addition of arrows};\hspace{5mm} $AB, \  BC  \ \mapsto \  AB+BC=AC$
	\vspace{2mm}
	
	\noindent Because of axiom A2 we can extend addition of arrows:
	
	\vspace{2mm}
	\noindent \textbf{generalized addition of arrows};\hspace{5mm} $AB+CD=AB+BX$, where $BX\sim CD$, under the condition that there is such a point $X$.
	\vspace{2mm}
	
	\textbf{Axiom} (\textbf{A3.1})  $AB\sim A'B'  \ \rightarrow \  BA\sim B'A'$.\\ 
\indent		(by the homogeneity principle)
	
	\textbf{Axiom} (\textbf{A3.2}) $AB\sim A'B' \ \land \ BC\sim B'C'  \ \rightarrow \  AC\sim A'C'$. \\
\indent		(by the homogeneity principle)
	
	Next axioms describe  multiplication of an arrow by a real number.
	
	\textbf{Axiom} (\textbf{A4})  $\forall \lambda , A, B \ \exists C  \ \ \lambda \cdot AB=AC$.\\
\indent		(by the very idea of the multiplication as stretching arrows) 
	
	\textbf{Axiom} (\textbf{A5})  $AB\sim CD \ \rightarrow \ \lambda AB\sim \lambda CD$.\\
\indent		(by the homogeneity principle) 
	
	\textbf{Axiom} (\textbf{A6.1})  $1\cdot AB = AB$.\\
\indent		(by the very idea of the multiplication as addition of the same arrow) 
	
	\textbf{Axiom} (\textbf{A6.2})  $ \lambda\cdot AB +\mu\cdot AB = (\lambda +\mu)\cdot AB  $.\\
\indent			(by the homogeneity principle and by the very idea of\\
\indent	 the multiplication as iterative addition of the same arrow)
	
	\textbf{Axiom} (\textbf{A6.3})  $\lambda\cdot (\mu\cdot AB)=(\lambda\cdot\mu)\cdot AB$. \\
\indent		(by the very idea of the multiplication \\ 
\indent	as iterative addition of the same arrow)
	
	\textbf{Axiom} (\textbf{A7}) (\textit{the scale invariance axiom})\\
	\indent  \textit{If} $ AC = \lambda\cdot AB$ \textit{and} $ AC' = \lambda\cdot AB'$ \textit{then} $CC' \sim \lambda\cdot BB'$.\\
\indent		(by the scale invariance principle)


\begin{figure}[h]
	\begin{center}
		\includegraphics[height=4cm]{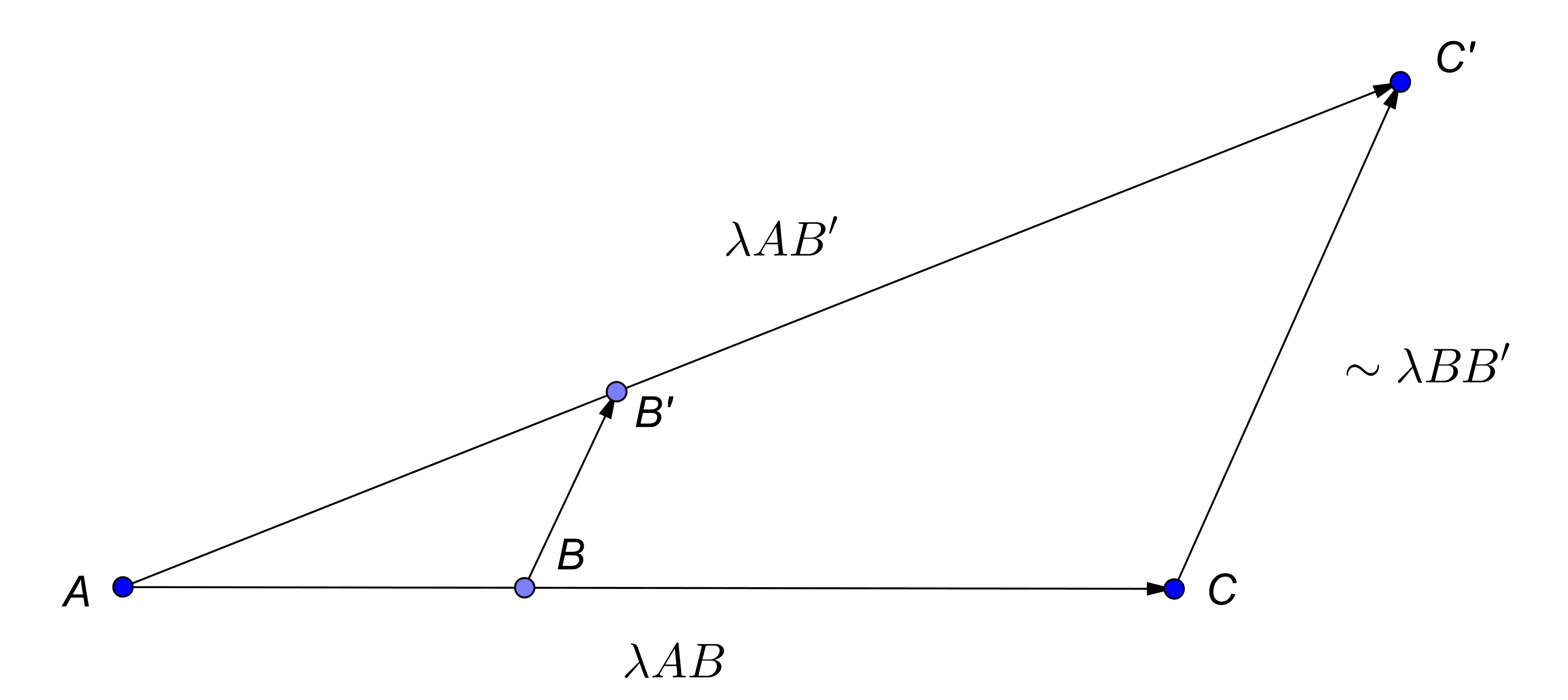}
		\caption{ }
		\label{fig:20}
	\end{center}
\end{figure}
	
	Next axioms describe  the distance function.
	
	\textbf{Axiom} (\textbf{A8})  $AB\sim CD \ \rightarrow\ |AB|=|CD|$.\\
\indent		(by the homogeneity principle)
	
	\textbf{Axiom} (\textbf{A9.1})  $|AA|=0$.\\ 
\indent		(by the very idea of measuring distance)
	
	\textbf{Axiom} (\textbf{A9.2})   $B\neq A  \ \rightarrow\ |AB|>0$. (\textit{positive definiteness})\\
\indent		(by the isotropy principle) 
	
	\textbf{Axiom} (\textbf{A9.3})   $|AB|=|BA|$.\\
\indent		(by the homogeneity principle) 
	
	\textbf{Axiom} (\textbf{A10})    $|\lambda AB|=\lambda |AB|$, \textit{for} $\lambda >0$. \\
\indent		(by the very idea of measuring along $AB$)
	
	\noindent The definition of a circle:
	
	\vspace{2mm}
	\noindent the \textbf{circle} with center $S$ and radius $r$ is $C(S,r)=\{ T: |ST| = r\}$
	\vspace{2mm}
	
	\textbf{Axiom} (\textbf{A11})    \textit{If a line has two   common points with a circle, points $A$ \\
	\indent	and $B$, then the midpoint $P$ of $AB$  is the point on the line nearest to the \\
	\indent center of the circle.}\\
\indent		(by the isotropy principle and an idea of continuity of space) 
	
	\textbf{Axiom} (\textbf{A12})   \textit{If a line has exactly one common point with a circle,\\
	\indent then the common point is the point on the line nearest to the center of \\
	\indent the circle. }\\
\indent		(by the isotropy principle and an idea of continuity of space)
	
	\noindent The definitions of perpendicularity of lines, orthogonal projection and scalar orthogonal projection of an arrow:
	
	\vspace{2mm}
	\noindent A line  $a$ is \textbf{perpendicular} to a line $b$, in symbols $a\bot b$, if $a$ intersects  $b$ and there is a point  $S$ on $a$ which is not on $b$ such that the intersection of $a$ and $b$ is the point on $b$ nearest to $S$.
	
	\vspace{2mm}
	\noindent\textbf{orthogonal projection} of a point $S$ to a line $p$ is the point on line $p$ nearest to the point $S$.
	\vspace{2mm}
	
	\noindent The \textbf{scalar orthogonal projection} of the arrow  $AB$ onto the arrow $CD$, denoted $AB_{CD}$, is the $\pm$ length of the orthogonal projection of the arrow $AB$ onto the line $p(CD)$,  where the sign is $+$ if the projection is in the direction of $CD$, $-$ otherwise.
	
	\textbf{Axiom} (\textbf{A13})   $|AB|=|AC| \ \rightarrow \ AB_{AC}=AC_{AB}$. \normalfont{(Fig.\ref{fig:20})} \\
\indent		(by the isotropy principle)
	
	\begin{figure}[h]
		\begin{center}
			\includegraphics[height=4cm]{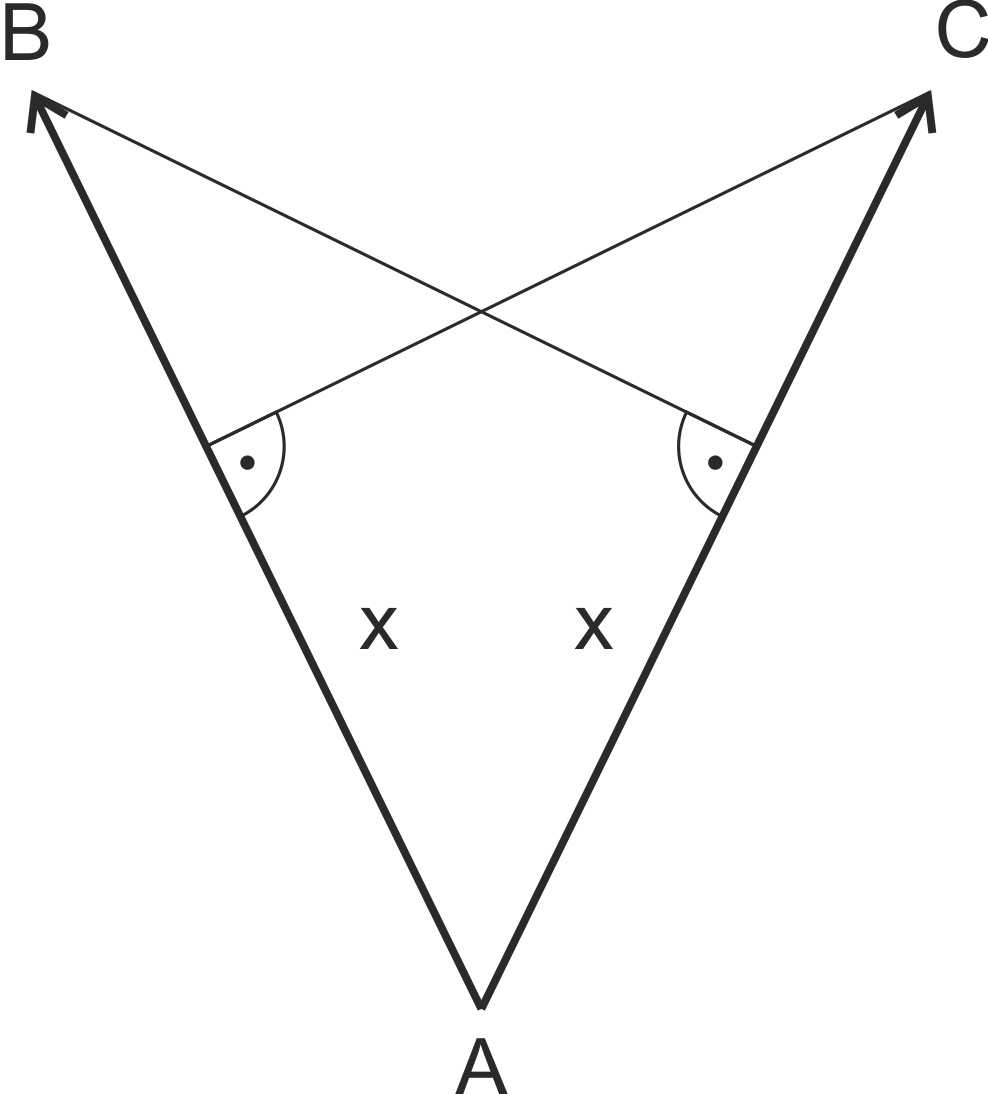}
			\caption{ }
			\label{fig:21}
		\end{center}
	\end{figure}

\bibliography{Euklid}
\bibliographystyle{alpha}


\end{document}